\documentclass[a4paper, 11pt]{amsart}
\usepackage{amsmath, amssymb, amsthm, amsfonts,enumitem,tikz,relsize}
\usepackage[normalem]{ulem}
\usepackage[colorlinks = true, citecolor = blue, hypertexnames=false]{hyperref}
\usepackage{xcolor}
\hypersetup{
    colorlinks,
    linkcolor={red!70!black},
    citecolor={blue!70!black},
    urlcolor={blue!80!black}
}
\usepackage{framed,booktabs,float,caption,scalerel,diagbox}
\usepackage[OT1]{fontenc}
\usepackage{comment}

\theoremstyle{plain}
\newcounter{MainTheoremCounter}
\newtheorem{MainTheorem}[MainTheoremCounter]{Theorem}

\newtheorem{theorem}{Theorem}[section]
\newtheorem{corollary}[theorem]{Corollary}
\newtheorem{lemma}[theorem]{Lemma}

\newtheorem{fact}[theorem]{Fact}
\newtheorem*{lemma*}{Lemma}
\newtheorem*{corollary*}{Corollary}

\theoremstyle{definition}
\newtheorem{definition}[theorem]{Definition}
\newtheorem{example}[theorem]{Example}

\newtheorem{remark}[theorem]{Remark}
\newtheorem{remarks}[theorem]{Remarks}

\newtheorem*{definition*}{Definition}
\newtheorem*{fact*}{Fact}
\newtheorem*{counterexample*}{Counterexample}

\newtheorem{task}{Remaining Task}
\newtheorem*{conjecture*}{Conjecture}
\newtheorem*{question*}{Question}
\newtheorem*{problem*}{Main Problem}

\newtheorem*{notation}{Notation}
\theoremstyle{remark}
\newtheorem{claim}{Claim}

\newenvironment{proofclaim}[1][Proof of Claim]{\begin{proof}[#1]}{\end{proof}}

\newcommand{\qoplus}{\mathbin{\scalerel*{\left(\mkern-2mu+\mkern-2mu\right)}{\bigcirc}}}

\newcommand{\cU}{\mathcal{U}}

\DeclareMathOperator{\Ar}{Ar}
\DeclareMathOperator{\Ob}{Ob}

\DeclareMathOperator{\ldim}{ldim}
\DeclareMathOperator{\ad}{ad}
\DeclareMathOperator{\tr}{tr}
\DeclareMathOperator{\im}{im}

\DeclareMathOperator{\charac}{char}
 
\DeclareMathOperator{\Aut}{Aut}

\DeclareMathOperator{\End}{End} 
\DeclareMathOperator{\Sym}{Sym} 
\DeclareMathOperator{\Alt}{Alt}
\DeclareMathOperator{\std}{std}

\DeclareMathOperator{\sgn}{sgn}
\DeclareMathOperator{\rstd}{r-std}
\DeclareMathOperator{\perm}{perm}

\DeclareMathOperator{\icosa}{Icosa}

\DeclareMathOperator{\GL}{GL}
\DeclareMathOperator{\PSL}{PSL} 
 
\DeclareMathOperator{\SL}{SL}

\DeclareMathOperator{\Nat}{Nat} 
\DeclareMathOperator{\Ad}{Ad}

\newcommand{\Mod}{\mathbf{Mod}}

\newcommand{\bF}{\mathbb{F}}

\newcommand{\bN}{\mathbb{N}}
\newcommand{\bZ}{\mathbb{Z}}

\usepackage[textsize=scriptsize, bordercolor=lightgray, linecolor=lightgray, backgroundcolor=white]{todonotes}

\begin{document}
\title[On $\Alt(n)$-modules when $n\le6$]{On $\Alt(n)$-modules with an additive dimension when $n\le6$}
\author{Barry Chin}
\email{\nolinkurl{barry_chin@live.com}}
\author{Adrien Deloro}
\address{Sorbonne Universit\'e and Universit\'e de Paris, \textsc{cnrs}, \textsc{imj-prg}, F-75005 Paris, France}
\email{adrien.deloro@imj-prg.fr}
\author{Joshua Wiscons}
\address{Department of Mathematics and Statistics\\
California State University, Sacramento\\
Sacramento, CA 95819, USA}
\email{joshua.wiscons@csus.edu}
\author{Andy Yu}
\address{Department of Mathematics and Statistics\\
California State University, Sacramento\\
Sacramento, CA 95819, USA}
\email{andyyu@csus.edu}

\date{\today}

\begin{abstract}
    Working in the general context of ``modules with an additive dimension,'' we complete the determination of the minimal dimension of a faithful $\Alt(n)$-module and classify those modules in three of the exceptional cases:  $2$-dimensional $\Alt(5)$-modules in characteristic $2$, $3$-dimensional $\Alt(5)$-modules in characteristic $5$, and $3$-dimensional $\Alt(6)$-modules in characteristic $3$. We also highlight the remaining work needed to complete the classification of the faithful $\Alt(n)$-modules of minimal dimension for all $n$; these open problems seem well suited as projects for advanced undergraduate or master's students. 
\end{abstract}

\maketitle

\section{Introduction}\label{S.Intro}
The setting of \emph{modules with an additive dimension} was recently introduced by Corredor and two of the present authors in an effort to unify and extend the existing theory of $G$-modules in various  contexts that support a well-behaved notion of dimension, \cite{CDW23}. For example, it treats simultaneously the classical setting of finite-dimensional vector spaces as well as certain ``tame'' model-theoretic classes of modules, such as those possessing finite Morley rank and those that are $o$-minimal. This setting provides a promising context in which to undertake  a  study of \emph{model-theoretic representation theory} with the work of \cite{CDW23}, and to some extent the present paper, establishing  basic results. Notably, this level of generality does not support character theory, nor basic tools such as Maschke's theorem, but ``elementary'' methods focusing on generators and relations still yield strong (and perhaps clarified) conclusions.

In \cite{CDW23}, the authors classify the faithful 
$\Sym(n)$- and $\Alt(n)$-modules of least dimension in this new context, assuming $n$ is large enough, hence solving the generic version of the following.

\begin{problem*}
Among modules with an additive dimension, identify the faithful $\Sym(n)$- and $\Alt(n)$-modules of least dimension for each $n$. 
\end{problem*}
However, \cite{CDW23} leaves open the topic for several small values of $n$; the left side of Table~\ref{tab.SmallDimension} provides a summary of the situation for $\Alt(n)$-modules. These remaining cases are complicated by exceptional isomorphisms between small instances of families of finite simple groups (e.g.~$\Alt(5) \cong \SL(2, 4)$). Via such isomorphisms, new geometric objects arise, and far from being exotic, the resulting modules are geometrically natural.
Our original focus for the present work was on $\Alt(n)$-modules with $n\ge 5$; however, our intermediate results also address smaller $n$ and the symmetric case.  

\begin{table}

\begin{minipage}{.49\linewidth}
{\small
\[\begin{array}{r|cccc} 
\hbox{\footnotesize \diagbox[width = 0.8cm,height = 0.8cm]{$n$}{$p$}}  
    & 2        & 3        & 5        & >5 \text{ or } 0 \\  \hline
\\[-10pt]
5   & -        & -        & -        & - \\[2pt] 
6   & -        & -        &  -       & - \\[2pt]
7   & -        & \fbox{s} & \fbox{s} & \fbox{s} \\[2pt]
8   & -        & \fbox{s} & \fbox{s} & \fbox{s} \\[2pt]
9   & 8        & \fbox{s} & \fbox{s} & \fbox{s} \\[2pt]
>9 & \fbox{s} & \fbox{s} & \fbox{s} & \fbox{s} \\
\end{array}\]
}
\caption*{Results of \cite{CDW23}}
\end{minipage}
\begin{minipage}{.49\linewidth}
{\small
\[\begin{array}{r|cccc}
\hbox{\footnotesize \diagbox[width = 0.8cm,height = 0.8cm]{$n$}{$p$}}   
    & 2        & 3        & 5        & >5 \text{ or } 0 \\  \hline
\\[-10pt]
5   & \fbox{2} & 3        & \fbox{s}        & 3 \\[2pt] 
6   & 4        & \fbox{3} &  5       & 5 \\[2pt]
7   & 4        & \fbox{s} & \fbox{s} & \fbox{s} \\[2pt]
8   & 4        & \fbox{s} & \fbox{s} & \fbox{s} \\[2pt]
9   & 8        & \fbox{s} & \fbox{s} & \fbox{s} \\[2pt]
>9 & \fbox{s} & \fbox{s} & \fbox{s} & \fbox{s} \\
\end{array}\]
}
\caption*{Incorporating our results}
\end{minipage}
\caption{Minimal dimensions of faithful $\Alt(n)$-modules in characteristic $p$: ``s'' stands for (reduced) ``standard'' (see \S\S~\ref{ss.FamiliarModules}) and as a dimension translates to $n-2$ if $p$ is a prime dividing $n$ and as $n-1$ otherwise. An entry in the table is boxed if the modules of least dimension have been classified.}
\label{tab.SmallDimension}
\end{table}
 
In this paper, we make  progress on the Main Problem by determining the value of the minimal dimension of a faithful $\Alt(n)$-module for all $n$ that are not treated in \cite{CDW23} and also provide identification of the minimal modules  when $(n,p)$ is $(5,2)$, $(5,5)$, or $(6,3)$. This is summarized on the right side of Table~\ref{tab.SmallDimension}; precise results and additional cases  are given in \S\S~\ref{ss.Results}.

A key take-away from our work echos that from \cite{CDW23}: even though the setting of modules with an additive dimension is quite general (with known objects that are not in the category of vector spaces), the minimal objects are indeed just the familiar ones.

\subsection{Examples}\label{ss.FamiliarModules}
Here we lay out some of the classical $\Sym(n)$- and $\Alt(n)$-modules that provide the desired lower bounds on the dimension of a faithful module and will be the targets for identification. We begin with the standard module; notation mostly coincides with that in \cite{CDW23}.

\begin{example}[Standard Module]\label{exam.standard}
Let $\perm(n,\bZ) = \bZ e_1\oplus \cdots \oplus \bZ e_n$ be the $\bZ[\Sym(n)]$-module where $\Sym(n)$ permutes the $e_i$ naturally (by permuting the subscripts). This is not a minimal $\Sym(n)$-module; two submodules are:
\begin{itemize}
\item 
$\std(n, \bZ) := \{\sum_i c_i e_i \mid c_i\in \bZ \text{ and } \sum c_i = 0\}$;
\item
$Z(\perm(n, \bZ)) := \{\sum_i ce_i \mid \text{$c\in \bZ$}\}$.
\end{itemize}

We can easily generalize these modules to have coefficients in any abelian group $L$. Indeed, view $L$ as a trivial $\Sym(n)$-module, and set $\perm(n,L) := \perm(n,\bZ) \otimes_\bZ L$, which has submodules:
\begin{itemize}
\item
$\std(n,L) := \std(n,\bZ) \otimes_\bZ L$;
\item
$Z(\perm(n, L)) := Z(\perm(n, \bZ)) \otimes_\bZ L$.
\end{itemize}

In the case when $Z(\perm(n, L)) \cap \std(n,L) \neq 0$, we can get a `smaller' module by considering the quotient $\rstd(n, L) = \std(n,L)/(Z(\perm(n, L)) \cap \std(n,L))$. We refer to $\std(n,L)$ as the \textbf{standard module} over $L$ and the quotient $\rstd(n, L)$ as the \textbf{reduced standard module} over $L$; these are of course modules for $\Alt(n)$ as well.
\end{example}

\begin{remarks}\hfill
\begin{enumerate}
    \item In the familiar linear setting when $L = F_+$ is the additive group of a field $F$, let $\ldim_F$ denote the linear dimension over $F$. Then $\ldim_F (\std(n,F_+)) = n-1$. Further, $\std(n,F_+)=\rstd(n, F_+)$ if and only if $\charac F \nmid n$; if $\charac F \mid n$, we have $\ldim_F (\rstd(n,F_+)) = n-2$. 
    \item In general, $\std(n,L)=\rstd(n, L)$ if and only if $L$ does \emph{not} have elements of order dividing $n$.
\end{enumerate}   
\end{remarks}

\begin{example}[Sign and Sign-tensored Modules]
Let $\sigma : \Sym(n) \rightarrow \{\pm 1\}$ denote the \textbf{sign} or \textbf{signature} homomorphism: $\sigma(g) = 1$ if and only if $g$ is an even permutation. We let $\sgn(n,\bZ) = \bZ$ be the $\Sym(n)$-module with action $g(v) = \sigma(g)\cdot v$; we refer to this as the \textbf{sign module}. Though we will not use it, we can ``extend the scalars'' to include any abelian group $L$ by letting $\sgn(n,L) = \sgn(n,\bZ) \otimes_\bZ L$.

The sign module is clearly not faithful when $n\ge 3$ (since $\Alt(n)$ acts trivially), but we can use it to build an important, and typically faithful, variant of the (reduced) standard module. For any abelian group $L$, we set 
\begin{itemize}
\item
$\std^\sigma(n,L) :=  \sgn(n,\bZ) \otimes_\bZ \std(n,L)$
\item
$\rstd^\sigma(n,L) :=  \sgn(n,\bZ) \otimes_\bZ \rstd(n,L)$.
\end{itemize}
We refer to these as the \textbf{sign-tensored (reduced) standard modules}. One effect of this construction is that an element of $\std(n,L)$ is centralized (i.e.~fixed) by a transposition $\tau$ if and only if the same element is inverted by $\tau$ in $\std^\sigma(n,L)$. Note that the difference between $\std(n,L)$ and $\std^\sigma(n,L)$ is only detectable by $\Sym(n)$; they are isomorphic as $\Alt(n)$-modules (as are $\rstd(n,L)$ and $\rstd^\sigma(n,L)$).
\end{example}   

The $\Sym(n)$- and $\Alt(n)$-modules of least dimension (in the classical setting and in the one  studied here) are always essentially of the form $\rstd(n,L)$ or $\rstd^\sigma(n,L)$ provided $n$ is large enough (see~\cite{CDW23}). But for small $n$, there are various exceptional modules, some of which we highlight now. 

\begin{example}[Exceptional case: $\Alt(5)$-modules of dimension $2$]\label{eg.NatSL2}
The main point is that $\Alt(5)\cong\SL(2,\bF_4)$. Now, for $F$ a field, we let $\Nat(n,F) = F^n$ be the \textbf{natural module} for $\SL(n,F)$ with elements acting by matrix multiplication. And, if $L$ is any vector space over $F$, we may, as before, extend $\Nat(n,F)$ to the $\SL(n,F)$-module $\Nat(n,L) =\Nat(n,F)\otimes_{F} L$. With this notation, we have that $\Nat(2,L)$ is an $\Alt(5)$-module whenever $L$ is a vector space over $\bF_4$.
\end{example}

\begin{example}[Exceptional case: $\Alt(6)$-modules of dimension $3$]\label{eg.AdPSL2}
Here we use that $\Alt(6)\cong\PSL(2,\bF_9)$. For $F$ a field, we let $\Ad(n,F)$ denote the set of $n\times n$ matrices over $F$ of trace $0$, viewed as a $\PSL(n,F)$-module with elements acting by conjugation; this is the \textbf{adjoint module} for $\PSL(n,F)$. As usual, we also define the module $\Ad(n,L) =\Ad(n,F)\otimes_{F} L$ for $L$ any vector space over $F$. So, in this notation, $\Ad(2,L)$ is an $\Alt(6)$-module whenever $L$ is a vector space over $\bF_9$.

One also has $\Alt(5)\cong\PSL(2,\bF_5)$, which gives rise to $3$-dimensional $\Alt(5)$-modules of the form $\Ad(2,L)$ for $L$ any vector space over $\bF_5$. These modules turn out to be isomorphic to $\rstd(5,L)$ (see Remark~\ref{r:actuallystandard} below).
\end{example}

\subsection{Main results}\label{ss.Results}
We  prepare to state our main results. As in \cite{CDW23}, we use $\Mod(G,d)$ to denote the collection of all \emph{dim-connected $G$-modules that carry an additive dimension}; a proper definition of being a $G$-module with an additive dimension is given in \S\S~\ref{s.ModAddDim} as is the definition of dim-connectedness. Both notions are quite natural, and in the classical setting of finite-dimensional vector spaces, all subspaces are automatically dim-connected. 
We say that a module $V$ has prime \emph{characteristic} $p$ if $V$ is an elementary abelian $p$-group; the notion of characteristic $0$ in our more general setting will be defined later (see Definition~\ref{d.Characteristic}). 

Our first theorem determines the bounds missing from the left side of Table~\ref{tab.SmallDimension} and also gives information about the case when $n= 4$. The completed Table~\ref{tab.SmallDimension} then  follows from combining it with \cite{CDW23} (see Remark~\ref{r.ProofOfTable}). 

\begin{MainTheorem}[Bounds]\label{t.AltBounds} 
Let $V\in \Mod(\Alt(n),d)$ be faithful.
\begin{enumerate}
\item If $n = 4$, then $d \ge 2$; if $d=2$, then $V$ has an  $\Alt(4)$-submodule of characteristic $2$ of positive dimension.
\item If $n = 5$, then $d \ge 2$; if $d=2$, then $V$ has characteristic $2$.
\item If $n = 6$, then $d \ge 3$. Moreover, if $d=3$, then $V$ has characteristic $3$, and if $d=4$, then $V$ has characteristic $2$ or $3$.
\end{enumerate}
\end{MainTheorem}

Our second theorem identifies the minimal $\Alt(n)$-modules in three special cases: those associated to the exceptional isomorphisms $\Alt(5)\cong\SL(2,\bF_4)$, $\Alt(5)\cong\PSL(2,\bF_5)$, and $\Alt(6)\cong\PSL(2,\bF_9)$; the target modules were defined in \S\S~\ref{ss.FamiliarModules}. We use $\Mod(G,d,p)$ to denote the subclass of $\Mod(G,d)$ consisting of those modules of characteristic $p$.

\begin{MainTheorem}[Recognition]\label{t.AltRec} 
\hfill

\begin{enumerate}
\item If $V\in \Mod(\Alt(5),2,2)$ is faithful, then $V\cong \Nat(2,L)$ for some $1$-dimensional subgroup $L\le V$.
\item If $V\in \Mod(\Alt(5),3,5)$ is faithful, then $V\cong \Ad(2,L)$ for some $1$-dimensional subgroup $L\le V$. 
\item If $V\in \Mod(\Alt(6),3,3)$ is faithful, then $V\cong \Ad(2,L)$ for some $1$-dimensional subgroup $L\le V$. 
\end{enumerate}
\end{MainTheorem}

\begin{remark}\label{r:actuallystandard}
Since  $\rstd(5,L)\in \Mod(\Alt(5),3,5)$, we indirectly see that $\Ad(2,L)\cong \rstd(5,L)$, so  modules in $\Mod(\Alt(5),3,5)$ are indeed standard as indicated in Table~\ref{tab.SmallDimension}.
\end{remark}

Though our main results focus on the alternating group, we also analyze the full symmetric group in the smallest of cases.  Lemmas~\ref{l.Sym3} and \ref{l.Sym3Recognition} identify the minimal faithful $\Sym(3)$-modules that are $3$-divisible; Lemmas~\ref{l.Sym4} and \ref{l.Sym4Recognition} do the same for  $\Sym(4)$-modules that are $2$-divisible. Recall that an abelian group $V$ is \emph{$p$-divisible} if for all $v\in V$ there is a $w \in W$ such that $v = pw$; this is a strong form of ``characteristic not $p$'' (see~Fact~\ref{fact.Divisibility}).

\subsection{Remaining work}\label{ss.RemainingWork}

A quick glance at Table~\ref{tab.SmallDimension} 
shows that, although the minimal dimensions of faithful $\Alt(n)$-modules are known, classification of the minimal modules 
remains open in certain cases. To solve the Main Problem one is left with four independent subtasks.

\begin{task}\label{p.Alt(5)}
Classify the minimal faithful $\Alt(5)$-modules in characteristics other than $2$ and $5$. 
\end{task}

This leads to the following question, which is a current favorite of ours: 
is every faithful module in $\Mod(\Alt(5),3,p)$ with $p\neq2,5$ derived from an ``icosahedron module'' of dimension $3$?
Let $\phi \in \mathbb{C}$ be the golden ratio.
Working in $(\bZ[\phi])^3$, we may build a regular icosahedron $I$ on the 12 vertices arising from all cyclic permutations of $(0,\pm1,\pm\phi)$. Viewing $\Alt(5)$ as the rotational symmetry group of $I$ then allows us to define an action of $\Alt(5)$ on $(\bZ[\phi])^3$; the resulting module we denote by $\icosa(5,\bZ[\phi])$.  We know of no general recognition theorem for $\icosa(5,\bZ[\phi])$ (analogous, say, to Fact~\ref{fact.Nat(SL2)}) and would be very happy to see one developed.
An interesting complication will be that not all finite fields have an analogue of $\phi$.
Do note that in the classical setting $\Mod(\Alt(5),3,p)$ with $p\neq2,5$ contains \emph{two} faithful modules: the second being derived from $\icosa(5,\bZ[\phi])$ via the Galois automorphism exchanging $\sqrt{5}$ and $-\sqrt{5}$.

\begin{task}\label{p.Alt(6)}
Classify the minimal faithful $\Alt(6)$-modules in characteristic not $3$.
\end{task}

Here we would first ask if every faithful module in $\Mod(\Alt(6),4,2)$ is of the form $\rstd(6,L)$ or $\rstd^\chi(6,L)$ where $\rstd^\chi(6,L)$ denotes the module obtained from $\rstd(6,L)$ via the outer automorphism $\chi\in \Aut(\Sym(6))$ applied to $\Alt(6)$. And for $p\neq2,3$, the question is if every faithful module in $\Mod(\Alt(6),5,p)$ is of the form $\std(6,L)$ or $\std^\chi(6,L)$. Actually, we should slightly modify this in characteristic $0$ to allow for a faithful $V\in \Mod(\Alt(6),5,0)$ to be \emph{dim-isogenous} to $\std(6,L)$ or $\std^\chi(6,L)$, which we take to mean that there is a surjective $\Alt(6)$-module homomorphism with a (possibly nontrivial) kernel of dimension $0$; this is typically what happens in \cite[Theorem]{CDW23}.

One strategy for approaching Task~\ref{p.Alt(6)} would be to revisit the Extension Lemma of \cite{CDW23} and show that some version of it holds for $n=6$ when the characteristic is not $3$. The main issue seems to be in adapting Claim~1 of the Extension Lemma (which makes significant use of $n\ge7$).

\begin{task}\label{p.AltInChar2}
Classify the minimal faithful $\Alt(n)$-modules in characteristic $2$ for $n=7,8,9$.
\end{task}

A key point for this problem is that $\Alt(8) \cong \SL(4,2)$; thus the faithful modules in $\Mod(\Alt(8),4,2)$ are expected to be of the form  $\Nat(4,L)$ or $\Nat^\psi(4,L)$ where $\Nat^\psi(4,L)$ denotes the module obtained from $\Nat(4,L)$ via the inverse-transpose (outer) automorphism of $\SL(4,2)$.  The case of $n=9$ requires separate treatment and presents an increased challenge: here there are at least three types of minimal modules~\cite{atlas3}.

\begin{task}\label{p.Sym(5)}
Classify the minimal faithful $\Sym(n)$-modules in characteristics not dividing $n$ for $n=5,6$.
\end{task}
This last task highlights the remaining gap in determining the minimal $\Sym(n)$-modules. This is raised at the beginning of \S\S~1.2 in \cite{CDW23} where known complications are also briefly mentioned.

\section{Background and preliminary results}

\subsection{Modules}\label{ss.modules}
Although the setting of modules with an additive dimension was introduced in \cite{CDW23}, many of our methods are classical. Here we collect a handful of basic tools for analyzing $G$-modules that will be used in our analysis; though typically quite classical, we include proofs to illustrate the types of arguments used later. For now, we do not assume that there is any notion of dimension around: a $G$-module is simply an abelian group $V$ and on which $G$ acts as automorphisms.

\begin{notation}
Let $G$ be a group. Assume $V$ is a $G$-module and $g\in G$. We view $g$ as an element of $\Aut(V)$, the automorphism group of $V$;  $\End(V)$ denotes the endomorphism ring of $V$. When $g$ has finite order $m$, we define maps $\ad_g$ and $\tr_g$ in $\End(V)$ as follows:
\begin{itemize}
\item
$\ad_g = 1-g$;
\item
$\tr_g = 1 + g + \dots + g^{m-1}$.
\end{itemize}
We refer to these as the \textbf{adjoint} and \textbf{trace} maps. They commute.
If the module $V$ is clear from the context, we use notation $B_g = \ad_g(V)$ and $C_g = \tr_g(V)$ for the images of $V$ under these maps.
When applying $g$ to an element $v\in V$, we typically write $gv$ in place of $g(v)$; however, we tend to include the parentheses when applying the adjoint and trace maps.
\end{notation}


\begin{lemma}\label{l.AdjointTrace}
Let $G$ be a group. Assume $V$ is a $G$-module and $g\in G$ has finite order. Then
\begin{enumerate}
    \item\label{l.AdjointTrace.i:BTrivial} for $X\subseteq V$, $\ad_g(X) = 0$ if and only if $g$ centralizes $X$;
    \item\label{l.AdjointTrace.i:CCentral} $\ad_g(C_g) = 0$, so $C_g$ is centralized by $g$;
    \item\label{l.AdjointTrace.i:traceonB} $\tr_g(B_g) = 0$, so if $|g| = 2$, then $B_g$ is inverted by $g$;
    \item\label{l.AdjointTrace.i:BCIntersection} every nontrivial element of $B_g \cap C_V(g)$ has order dividing $|g|$;
    \item\label{l.AdjointTrace.i:AdTrConjugation} for all $h\in G$, $h\circ \ad_g\circ h^{-1} = \ad_{hgh^{-1}}$ and $h\circ \tr_g\circ h^{-1} = \tr_{hgh^{-1}}$. 
\end{enumerate}
\end{lemma}
\begin{proof}
Set $|g|=m$. We proceed point by point. 
\begin{enumerate}
\item[\eqref{l.AdjointTrace.i:BTrivial}] This follows  from the definition of the adjoint map: $\ad_g(X) = 0$ if and only if for all $v\in X$ one has $v-gv=0$ (equivalently $v= gv$).

\item[\eqref{l.AdjointTrace.i:CCentral}] For each $v\in C_{g}$ there is $w\in V$ such that $v=\left(1+g+\cdots+g^{m-1}\right) w$, so:
\[ \ad_g(v) 
= \left(1-g\right)\left(1+g+\cdots+g^{m-1}\right) w
=\left(1-g^m\right) w
=0  w.\]
Thus, $\ad_g(C_g) = 0$, which is equivalent to $g$ centralizing $C_g$ by  \eqref{l.AdjointTrace.i:BTrivial}.

\item[\eqref{l.AdjointTrace.i:traceonB}] Similar to the proof of \eqref{l.AdjointTrace.i:CCentral}, we may write each $v\in B_{g}$ as $v=\left(1-g\right) w$
for some $w\in V$, so
\[ \tr_g(v) 
= \left(1+g+\cdots+g^{m-1}\right)\left(1-g\right) w
=\left(1-g^m\right) w
=0  w.\]
Thus,  $\tr_g(B_g) = 0$. Further, if $|g| = 2$, then $0 = \tr_g(v) = v+g\cdot v$ if and only if $g\cdot v = -v$.

\item[\eqref{l.AdjointTrace.i:BCIntersection}] Let $v \in B_g \cap C_V(g)$. Since $v \in B_g$, one has $\tr_g(v) = 0$ by \eqref{l.AdjointTrace.i:traceonB}, meaning $(1 + g + \cdots + g^{m-1})v = 0$. Also, $v\in C_V(g)$, so $(1 + g + \cdots + g^{m-1})v = v + gv + \cdots g^{m-1}v = mv$. Thus, $mv = 0$, which implies $v$ has order dividing $m$. 

\item[\eqref{l.AdjointTrace.i:AdTrConjugation}] Observe that  $h\circ\ad_g\circ h^{-1} = h(1-g)h^{-1} = 1-hgh^{-1} =\ad_{hgh^{-1}}$, and similarly, 
\begin{align*}
 h\circ\tr_g\circ h^{-1} 
    & =h(1+g+\cdots+g^{m-1})h^{-1}\\
    & =(1+hgh^{-1}+\cdots+hg^{m-1}h^{-1})\\
    & =\left(1+hgh^{-1}+\cdots+\left(hgh^{-1}\right)^{m-1}\right)=\tr_{hgh^{-1}}.\qedhere
 \end{align*}
\end{enumerate}
\end{proof}

\begin{lemma}\label{l.PermutingBandC}
Let $G$ be a group. Assume $V$ is a $G$-module and $g\in G$ has finite order. Then
\begin{enumerate}
    \item\label{i:BCInversion} $B_g = B_{g^{-1}}$ and $C_g = C_{g^{-1}}$;
    \item\label{i:BCConjugation} for all $h\in G$, $h(B_g) = B_{hgh^{-1}}$ and $h(C_g) = C_{hgh^{-1}}$.
\end{enumerate}
In particular, if $h$ centralizes or inverts $g$, then $B_g$ and $C_g$  are $\langle h\rangle $-invariant.
\end{lemma}
\begin{proof}
Set $|g|=m$.
\begin{enumerate}
\item[\eqref{i:BCInversion}] As $g$ is an automorphism of $V$, $g(V) = V$, so also $-g(V) = g(-V) = g(V) = V$. Thus, 
$B_{g^{-1}} = (1-g^{-1})(V) = (1-g^{-1})(-g)(V) = (-g+1)(V) = B_g$. Similarly, $g^{m-1}$ is an automorphism of $V$, so $g^{m-1}(V) = V$. Thus, $C_{g^{-1}} = (1+g^{-1}+\cdots+g^{-(m-1)})(V) = (1+g^{-1}+\cdots+g^{-(m-1)})g^{m-1}(V)  = (g^{m-1}+g^{m-2}+\cdots+1)(V) = C_{g}$. 
\item[\eqref{i:BCConjugation}] Using Lemma~\ref{l.AdjointTrace}\eqref{l.AdjointTrace.i:AdTrConjugation} and the fact that $V = h^{-1}V$ since $h$ is an automorphism, we see that
\[
hB_{g} = h\circ \ad_g(V) = h\circ \ad_g\circ h^{-1}(V) = \ad_{hgh^{-1}}(V) = B_{hgh^{-1}},
\]
and a similar calculation shows $h(C_g) = C_{hgh^{-1}}$.
\end{enumerate}

For the final remark, if $h$ centralizes or inverts $g$, then $hB_{g}=B_{hgh^{-1}} = B_{g^{\pm1}} = B_g$ and similarly for $C_{g}$. So in this case, $B_g$ and $C_g$  are $\langle h\rangle $-invariant.
\end{proof}

\begin{notation}
Let $G$ be a group acting on an abelian group $V$. We define \[[\underbrace{G,\ldots,G}_n,V] = \bigoplus_{g_1,\ldots,g_n\in G} \ad_{g_1}\circ \cdots \circ \ad_{g_n} (V).\]
This agrees with the usual definition of the commutator  $[G,\ldots,G,V]$ (with appropriate bracketing) in the corresponding group $G\ltimes V$.
\end{notation}

\begin{lemma}\label{l.primality}
Let $V$ be a nontrivial elementary abelian $p$-group. If 
 $V$ is a $G$-module for $G$ a finite $p$-group, then $[G,V] < V$.
\end{lemma}
\begin{proof}
The main point is that if $|g|=p$, then $\ad_g(V) <V$. To see this, suppose $\ad_g(V) = V$. Then $(1-g)V = V$, so $(1-g)^pV = V$. But, the characteristic is $p$, so $V = (1-g)^pV = (1-g^p)V = 0$, a contradiction.

We now proceed by induction on $|G|$. If $|G|=p$, then  writing $G = \langle g \rangle$, we have that $[G,V] = [\langle g \rangle,V] = \ad_g(V) <V$.  Now suppose $|G|>p$; by nilpotence of finite $p$-groups (a phenomenon lost for infinite $p$-groups), $G$ has some proper nontrivial normal subgroup $N$. By induction, we have $[N,V] < V$. Since $N$ acts trivially on $V/[N,V]$,  $G$ acts on  $\overline{V} = V/[N,V]$ as $G/N$. Again, by induction, we have $[G/N,\overline{V}] < \overline{V}$, so $[G,\overline{V}] < \overline{V}$. This final point implies $[G,V] < V$.
\end{proof}

We will also require a fundamental fact about so-called quadratic actions. 

\begin{definition}
Let $U$ be a group acting on an abelian group $V$.
We say that the action is \textbf{quadratic} if $[U,U,V] = 0$ but  $[U,V] \neq 0$. 
\end{definition}

\begin{lemma}\label{l.Quadratic}
If a group $U$ acts faithfully and quadratically on an abelian group $V$, then $U$ is abelian. 
\end{lemma}
\begin{proof}
    Let $g,h\in U$ be arbitrary. By assumption,  $0 = \ad_g(\ad_h(v)) = \ad_g(v - h v) =    v-g v - h v + gh v$ for all $v\in V$. Since this holds for all $g$ and $h$, we may switch the order to get $0 = v - h v - g v +hg v$ as well. Together, these  equations imply that $gh v = 0 = hg v$, so as the action is faithful, it must be that $gh = hg$, as desired.
\end{proof}

\subsection{Modules with an additive dimension}\label{s.ModAddDim}
We now introduce modules with an additive dimension; our definitions are taken directly from \cite{CDW23}. We refer  to \cite[Section 2]{CDW23} for more details and strongly encourage the reader to review the motivation and various examples found there. Throughout, we use $\ker f$ and $\im f$ for the kernel and image of a map $f$.

\begin{definition}
A \textbf{modular universe} $\cU$ is a collection of abelian groups $\Ob(\cU)$ and homomorphisms $\Ar(\cU)$ between them  that satisfy the following closure properties.
\begin{itemize}
\item{}[\textsc{inverses}]
If $f \in \Ar(\cU)$ is an isomorphism, then $f^{-1} \in \Ar(\cU)$.
\item{}[\textsc{products}] If $V_1, V_2 \in \Ob(\cU)$ and $f_1, f_2 \in \Ar(\cU)$, then 
$V_1 \times V_2 \in \Ob(\cU)$, and 
	$\Ar(\cU)$ contains $f_1 \times f_2$,  the projections $\pi_i : V_1\times V_2 \rightarrow V_i$, and the diagonal embeddings $\Delta_k:V_1\rightarrow V_1^k$.
\item{}[\textsc{sections}] 
If $W \leq V$ are in $\Ob(\cU)$, then $V/W\in \Ob(\cU)$ and $\Ar(\cU)$ contains the inclusion $\iota:W\rightarrow V$ and quotient $p:V\rightarrow V/W$ maps.
\item{}[\textsc{kernels/images}] 
If $f:V_1 \rightarrow V_2$ is in $\Ar(\cU)$, then $\ker f, \im f \in \Ob(\cU)$, and for all $W_1,W_2 \in \Ob(\cU)$,
	\begin{itemize} 
	\item
	if $W_1 \leq \ker f$, the induced map $\overline{f}\colon V_1/W_1 \to V_2$ is in $\Ar(\cU)$;
	\item 
	if $\im f \le W_2 \le V_2$,  the induced map $\check{f}\colon V_1 \to W_2$ is in $\Ar(\cU)$.
	\end{itemize} 
\item{}[$\bZ$-\textsc{module structure}]
If $V \in \Ob(\cU)$, then $\Ar(\cU)$ contains the sum  map $\sigma \colon V\times V \to V$ and the multiplication-by-$n$  maps $\mu_n\colon V \to V$.
\end{itemize}
The groups in $\Ob(\cU)$ are called the \textbf{modules} of $\cU$ and the homomorphisms in $\Ar(\cU)$ its \textbf{compatible morphisms}.

If $V$ is a module in  $\cU$ and a group $G$ acts on $V$ by \emph{compatible} morphisms of $\cU$, we say that $V$ is a \textbf{$G$-module} in $\cU$.
\end{definition}

\begin{definition}\label{d.ModAddDim}
Let $\cU$ be a modular universe.
An \textbf{additive dimension} on $\cU$ is a function $\dim\colon \Ob(\cU) \to \bN$ such that for all $f\colon V \to W$ in $\Ar(\cU)$, 
\[\dim V = \dim \ker f + \dim \im f.\]

If $(\cU, \dim)$ is a modular universe with an additive dimension and $V \in \Ob(\cU)$, we call $V$ \textbf{a module with an additive dimension}, leaving $\cU$ implicit from context. 
\end{definition}

\begin{definition}
A module $V$ with an additive dimension is called \textbf{dim-connected} (for dimension-connected) if every proper submodule $W<V$ in $\cU$ satisfies $\dim W < \dim V$.
\end{definition}

\begin{fact}[{\cite[Connectedness Properties]{CDW23}}]\label{fact.ConnectednessProps}
Let $(\cU, \dim)$ be a modular universe with an additive dimension, and let $V_1,V_2$ be modules in $\cU$.
\begin{enumerate}
\item
If $V_1$ is dim-connected and $f\colon V_1 \to V_2$ is compatible, then $\im f$ is dim-connected.
\item
If $V_1$ and $V_2$ are dim-connected, then so is $V_1 \times V_2$.
\item
If $V_1, V_2 \leq V$ and $V_1$ and $V_2$ are dim-connected, then so is $V_1+V_2$.
\end{enumerate}
\end{fact}

\begin{notation}
Let $G$ be a group. We use $\Mod(G)$ to denote the class of all \emph{dim-connected} $G$-modules with an additive dimension (each coming from its own dimensional universe); $\Mod(G, d)$ denotes the subclass of those modules of dimension exactly $d$.
\end{notation}

\begin{fact}[{\cite[Divisibility Properties]{CDW23}}]\label{fact.Divisibility}
Let $V\in \Mod(G)$ for some group $G$. If $p$ is a prime, then $V$ is $p$-divisible (as defined in \S~\ref{S.Intro}) if and only if $\Omega_p(V) := \{v \in V: pv = 0\}$ has dimension $0$. 
\end{fact}

\begin{definition}
If $A_1,\ldots,A_n$ are submodules of a module $V$, the sum $\sum A_i$ is said to be \textbf{quasi-direct} if $\dim \sum A_i = \sum \dim A_i$, in which case we write $\sum A_i = A_1\qoplus \cdots\qoplus A_n$.
\end{definition}


\begin{fact}[{\cite[Coprimality Lemma]{CDW23}}]\label{fact.Coprimality}
Let $V\in \Mod(G)$ for a group $G$. Suppose $g\in G$ has prime order $p$ and $V$ is $p$-divisible. Then $V = B_g \qoplus C_g$.
\end{fact}

\begin{remark}\label{rem.dcComponent}
It follows from Fact~\ref{fact.Coprimality} that, under the same assumptions, $\dim C_V(g) = \dim C_g$; this is discussed in the Remarks following the proof of \cite[Coprimality Lemma]{CDW23}. Since $\dim C_V(g)$ has a 
dim-connected submodule of equal dimension, it has only one, and we say $C_g$ is the \emph{dc-component} of $C_V(g)$, see \cite[Remarks, page 11]{CDW23}.
\end{remark}

The next lemma presents the so-called Weight Lemma from \cite{CDW23}. We include a proof here to illustrate some of the methods we use later; it also presents a slightly different point of view than in \cite{CDW23}. 

\begin{lemma}[{cf.~\cite[Weight Lemma]{CDW23}}]\label{l.WeightLemma}
Let $G$ be a group with $K\trianglelefteq G$ a normal Klein $4$-subgroup, and let $V \in \Mod(G)$. Write $K=\{1,\alpha_1, \alpha_2, \alpha_3\}$. Define $V_{\alpha_1} = \tr_{\alpha_1} \circ \ad_{\alpha_2} \circ \ad_{\alpha_3}(V)$, $V_{\alpha_2} = \ad_{\alpha_1} \circ \tr_{\alpha_2} \circ \ad_{\alpha_3}(V)$, and $V_{\alpha_3} = \ad_{\alpha_1} \circ \ad_{\alpha_2} \circ \tr_{\alpha_3}(V)$;  set $C_K = \tr_{\alpha_1} \circ \tr_{\alpha_2} \circ \tr_{\alpha_3}(V)$.

If $V$ is $2$-divisible, then $V = V_{\alpha_1} \qoplus V_{\alpha_2} \qoplus V_{\alpha_3} \qoplus C_K$, and if the nontrivial elements of $K$ are conjugate in $G$, then $\dim V = 3\cdot\dim V_{\alpha_1} + \dim C_K$.
\end{lemma}
\begin{proof}
We have $V_{\alpha_{1}}=(1+\alpha_{1})(1-\alpha_{2})(1-\alpha_{3})V$ (similarly for $V_{\alpha_{2}}$ and $V_{\alpha_{3}}$) and $C_{K}=(1+\alpha_{1})(1+\alpha_{2})(1+\alpha_{3})V$.  Note that,
since the $\alpha_{i}$ commute for all $i$, the order of the factors in each expression does not matter.

Bear in mind the dim-connectedness of $B_{\alpha_i}$ and $C_{\alpha_i}$.

Suppose $V$ is $2$-divisible. Since $|\alpha_{i}|=2$, we may apply Fact~\ref{fact.Coprimality} to see that $V=B_{\alpha_{i}}\qoplus C_{\alpha_{i}}$.
Since $B_{\alpha_{i}}$ is a dim-connected submodule of $V$, it is
also $2$-divisible by Fact~\ref{fact.Divisibility}, so by Fact~\ref{fact.Coprimality} applied with $B_{\alpha_{i}}$ in place of $V$, $B_{\alpha_{i}}=\ad_{\alpha_{j}}(B_{\alpha_{i}})\qoplus\tr_{\alpha_{j}}(B_{\alpha_{i}})$. Similarly, $C_{\alpha_{i}}=\ad_{\alpha_{j}}(C_{\alpha_{i}})\qoplus\tr_{\alpha_{j}}(C_{\alpha_{i}})$. Thus, we find:
\begin{align*}
V &=B_{\alpha_{3}}(+)C_{\alpha_{3}}\\
&=\ad_{\alpha_{2}}(B_{\alpha_{3}})\qoplus \tr_{\alpha_{2}}(B_{\alpha_{3}})\qoplus \ad_{\alpha_{2}}(C_{\alpha_{3}})\qoplus \tr_{\alpha_{2}}(C_{\alpha_{3}})\\
&=\ad_{\alpha_{2}}\circ\ad_{\alpha_{3}}(V)\qoplus \tr_{\alpha_{2}}\circ\ad_{\alpha_{3}}(V)\qoplus \ad_{\alpha_{2}}\circ\tr_{\alpha_{3}}(V)\qoplus \tr_{\alpha_{2}}\circ\tr_{\alpha_{3}}(V).
\end{align*}

 We look at how $\alpha_{1}$ acts on the above four factors. Remember that $\alpha_{1}=\alpha_{2}\alpha_{3}$. Consider the first factor $W=\ad_{\alpha_{2}}\circ\ad_{\alpha_{3}}(V)$. Notice that $W\le B_{\alpha_{2}}\cap B_{\alpha_{3}}$, so
both $\alpha_{2}$ and $\alpha_{3}$ invert $W$ by Lemma~\ref{l.AdjointTrace}\eqref{l.AdjointTrace.i:traceonB}. Then $\alpha_{1}=\alpha_{2}\alpha_{3}$ must centralize $W$, so for all $w\in W$, $\tr_{\alpha_{1}}(w)=w+\alpha_{1}w=w+w=2w$.
Thus $\tr_{\alpha_{1}}(W)=2W$, so by 2-divisibility, $\tr_{\alpha_{1}}(W)=W$. This shows that
\begin{itemize}
    \item  $\ad_{\alpha_{2}}\circ\ad_{\alpha_{3}}(V) =\tr_{\alpha_{1}}\circ\ad_{\alpha_{2}}\circ\ad_{\alpha_{3}}(V)=V_{\alpha_{1}}.$
\end{itemize}
Similar calculations show
\begin{itemize}
    \item $\tr_{\alpha_{2}}\circ\ad_{\alpha_{3}}(V) = \ad_{\alpha_{1}}\circ\tr_{\alpha_{2}}\circ\ad_{\alpha_{3}}(V)=V_{\alpha_{2}}$;
    \item $\ad_{\alpha_{2}}\circ\tr_{\alpha_{3}}(V)=  \ad_{\alpha_{1}}\circ\ad_{\alpha_{2}}\circ\tr_{\alpha_{3}}(V)=V_{\alpha_{3}}$;
    \item $\tr_{\alpha_{2}}\circ\tr_{\alpha_{3}}(V) = \tr_{\alpha_{1}}\circ\tr_{\alpha_{2}}\circ\tr_{\alpha_{3}}(V)=C_{K}$.
\end{itemize}
Thus,
$
V=V_{\alpha_{1}}\qoplus V_{\alpha_{2}}\qoplus V_{\alpha_{3}}\qoplus C_{K}.
$

Now suppose that $\alpha_{1},\alpha_{2},\alpha_{3}$ are conjugate in $G$. For each  pair $i\neq j$ there exists $g\in G$ such that $g\alpha_{i}g^{-1}=\alpha_{j}$, and for $k\notin\{i,j\}$, $\{g\alpha_{j}g^{-1},g\alpha_{k}g^{-1}\} = \{\alpha_i,\alpha_k\}$ as $K$ is normal.
Then, similar to the proof of Lemma~\ref{l.PermutingBandC}, 
\begin{align*}
gV_{\alpha_{i}} & =g\circ\tr_{\alpha_{i}}\circ\ad_{\alpha_{j}}\circ\ad_{\alpha_{k}}(V)\\
& =g\circ\tr_{\alpha_{i}}\circ\ad_{\alpha_{j}}\circ\ad_{\alpha_{k}}\circ g^{-1}(V)\\
& =\tr_{g\alpha_{i}g^{-1}}\circ\ad_{g\alpha_{j}g^{-1}}\circ\ad_{g\alpha_{k}g^{-1}}(V)\\
& =\tr_{\alpha_{j}}\circ\ad_{g\alpha_{j}g^{-1}}\circ\ad_{g\alpha_{k}g^{-1}}(V)\\
 & =\tr_{\alpha_{j}}\circ\ad_{\alpha_{i}}\circ\ad_{\alpha_{k}}(V)=V_{\alpha_{j}}.
\end{align*}
Thus $g(V_{\alpha_{i}})=V_{\alpha_{j}}$, so $\dim V_{\alpha_{i}} = \dim V_{\alpha_{j}}$.  Since $V=V_{\alpha_{1}}\qoplus V_{\alpha_{2}}\qoplus V_{\alpha_{3}}\qoplus C_{K}$, 
$\dim V=3\cdot\dim V_{\alpha_{i}}+\dim C_{K}$.
\end{proof}

\begin{definition}\label{d.Characteristic}
Let $V$ be an abelian group. We say
\begin{itemize}
\item
 $V$ has \textbf{characteristic $p$}, for $p$ a prime, if $V$ is an elementary abelian $p$-group;
\item
$V$ has \textbf{characteristic $0$} if it is divisible, i.e.~if $V$ is $p$-divisible  for all primes $p$.
\end{itemize}
\end{definition}

\begin{notation}
Extending our earlier notation, we use $\Mod(G, d, p)$ to denote the subclass of $\Mod(G, d)$ consisting of those modules of characteristic $p$ (possibly equal to $0$).
\end{notation}

Not all modules have a well-defined characteristic, but it turns out that the ``irreducible'' ones do.

\begin{definition}\label{d.Irreducible}
Let $V\in \Mod(G)$ for some group $G$. We say that $V$ is \textbf{dc-irreducible} (for ``dimension-connected-irreducible''), if it has no non-trivial, proper, dim-connected $G$-submodule.
\end{definition}

\begin{fact}[{see \cite[Characteristic Lemma]{CDW23}}]\label{fact.Characteristic}
Let $V\in \Mod(G)$ for some group $G$. If $V$ is dc-irreducible, then $V$ has a characteristic.
\end{fact}

\section{Bounding minimal dimensions}

The bulk of our work occurs in this section. The main outcome will be a proof of Theorem~\ref{t.AltBounds}.

\begin{lemma}\label{l.Sym2}
Let $V\in \Mod(G,1)$ for some group $G$. If $g \in G$ has order $2$, then $g$ either centralizes or inverts $V$.
\end{lemma}
\begin{proof}
Consider $B_g$; it has dimension $0$ or $1$. If $\dim B_g = 0$, then as $B_g$ is dim-connected, $B_g = 0$, and  Lemma~\ref{l.AdjointTrace}\eqref{l.AdjointTrace.i:BTrivial} implies that $g$ centralizes $V$.  Now suppose that $\dim B_g= 1$; then dimension connectedness of $V$ implies $B_g=V$. Thus, in this case, Lemma~\ref{l.AdjointTrace}\eqref{l.AdjointTrace.i:traceonB} implies that $g$ inverts $V$.
\end{proof}

\begin{remark}\label{r.dimension1}
Lemma~\ref{l.Sym2} states that an element of order 2 must fix or invert a dim-connected, $1$-dimensional module. This is of course true for any element of order $2$, so if  $V\in \Mod(G,1)$ and $a,b\in G$ both have order $2$, then at least one of $a$, $b$, or $ab$ must centralize $V$.
\end{remark}

The next lemma (Lemma~\ref{l.Sym3}) as well as Lemma~\ref{l.Sym4} appear as  Claim~1 of \cite[First Geometrisation Lemma]{CDW23}, but we include proofs here for completeness and to add further detail. 

\begin{lemma}\label{l.Sym3}
If $V \in \Mod(\Sym(3), d)$ is faithful, then $d \ge 2$. 
\end{lemma}
\begin{proof}
Since $V$ is dim-connected, $d\ge 1$ as $d=0$ would imply $V=\{0\}$, contradicting faithfulness. Now suppose $d=1$. Let $\tau_1$ and $\tau_2$ be distinct transpositions in $\Sym(3)$. By Remark~\ref{r.dimension1}, one of $\tau_1$, $\tau_2$, or $\tau_1\tau_2$ fixes $V$. This contradicts faithfulness, so $d \ge 2$.
\end{proof}

\begin{corollary}\label{c.Bgamma}
Let $V \in \Mod(\Sym(3))$. If $\gamma \in \Sym(3)$ has order $3$, then any dim-connected, $1$-dimensional $\Sym(3)$-submodule of $\ad_\gamma(V)$ has characteristic $3$. In particular, if $V$ is $3$-divisible, then $\dim \ad_\gamma(V) \neq 1$.
\end{corollary}
\begin{proof}
Let $U$ be a dim-connected, $1$-dimensional $\Sym(3)$-submodule of $B_\gamma = \ad_\gamma(V)$.
By Lemma~\ref{l.Sym3}, the action of $\Sym(3)$ on $U$ is not faithful, so $\gamma$ must be in the kernel. In other words, $\gamma$ centralizes  $U$, so $U \le  C_V(\gamma)$. Thus,  $U\le B_\gamma \cap C_V(\gamma)$, which implies $U$ has characteristic $3$ by Lemma~\ref{l.AdjointTrace}\eqref{l.AdjointTrace.i:BCIntersection}.
 
Now, if $V$ is $3$-divisible, then $V$ has no subgroup of characteristic  $3$ of positive dimension by Fact~\ref{fact.Divisibility}.
In particular, cannot $\dim B_\gamma = 1$ happen.
\end{proof}

\begin{lemma}\label{l.Alt4}
If $V \in \Mod(\Alt(4), d)$ is faithful, then $d \ge 2$;  if $V$ is also $2$-divisible, then $d \ge 3$. 
\end{lemma}
\begin{proof}
Let $K = \left\{1, \alpha_{1},\alpha_{2},\alpha_{3}\right\}$ be the (normal) Klein 4-group in $\Alt(4)$. Then each $\alpha_i$ has order $2$, and $\alpha_{3} = \alpha_{1}\alpha_{2}$. If $d=1$, then by Remark~\ref{r.dimension1}, some $\alpha_i$ fixes $V$, so $V$ is not faithful. Thus $d\ge 2$.

Now suppose $d=2$ and, towards a contradiction, also assume $V$ is $2$-divisible. The three nontrivial elements of $K$ are conjugate in $\Alt(4)$, so
by Lemma~\ref{l.WeightLemma}, $\dim V=3\cdot\dim V_{\alpha_{1}}+\dim C_{K}$. Since the dimension of $V$ is $2$, the dimension of each  $V_{\alpha_{i}}$ is $0$, so  $\dim V = \dim C_{K}$. As $V$ is dim-connected, $V=C_K \le C_V(K)$, which contradicts faithfulness.
\end{proof}

\begin{remark}
Despite some effort, it remains unclear to us if a faithful $V \in \Mod(\Alt(4), 2)$  necessarily has characteristic $2$. Lemma~\ref{l.Alt4} (together with Fact~\ref{fact.Divisibility}) shows that such a $V$ must have exponent dividing $4$, but exponent equal to $4$ has yet to admit it does not exist. Of course, exponent $2$ does exist as can be seen via $\Alt(4) < \Alt(5) \cong \SL(2,\bF_4)$.

We remark that in \cite{CDW23} it is erroneously stated (in the Examples following the Characteristic Lemma)   that $V = (\mathbb{Z}/4\mathbb{Z})^2$ can be made a faithful $2$-dimensional $\Alt(4)$-module with action given by any embedding of $\Alt(4)$ into $\GL\left(2,\mathbb{Z}/4\mathbb{Z}\right)$. Certainly $V$ can be made a faithful $\Alt(4)$-module, but, as we show, it cannot have dimension $2$ (in the sense of Definition~\ref{d.ModAddDim}). 

The group $\GL\left(2,\mathbb{Z}/4\mathbb{Z}\right) \cong (\Alt(4) \rtimes\mathbb{Z}/4\mathbb{Z})\rtimes\mathbb{Z}/2\mathbb{Z}$ (as can be seen, for example, using GAP~\cite{GAP4}) contains a unique copy of $\Alt(4)$, which is given by   
\[H=\left\langle 
\left(\begin{array}{cc}
1 & 2\\
2 & 1
\end{array}\right),
\left(\begin{array}{cc}
2 & 1\\
1 & 1
\end{array}\right)\right\rangle.\]
Let $\alpha$ denote the first generator, and notice that $0 < C_\alpha < B_\alpha < V$. Thus, if this $V$ is in $\Mod(\Alt(4), d)$, we must have $d>2$.
\end{remark}

\begin{lemma}\label{l.Sym4}
If $V \in \Mod(\Sym(4), d)$ is faithful, then $d \ge 3$. 
\end{lemma}
\begin{proof}
Since $\Sym(4)$ contains a copy of $\Sym(3)$, Lemma~\ref{l.Sym3} implies that $d\ge2$. Similarly, Lemma~\ref{l.Alt4} can be used when $V$ is $2$-divisible to show that $d \ge 3$. 
Thus it remains to show that if $V$ is not $2$-divisible, then $d \neq 2$.

Assume that $d = 2$ and that $V$ is not $2$-divisible. Let $W=2V$. Then $W$ is a dim-connected submodule of $V$. Since $V$ is not $2$-divisible,  $W\neq V$. Thus, as $V$ is dim-connected, $\dim W \le 1$. 

If $\dim W = 1$, consider the action of $\Sym(4)$ on $W$ and $V/W$. By Remark~\ref{r.dimension1} applied to the transpositions, the product of any two transpositions fixes $W$, so $\Alt(4)$ acts trivially on $W$. Similarly, $\Alt(4)$ acts trivially on $V/W$. This implies that  $\Alt(4)$ acts quadratically on $V$, so by Lemma~\ref{l.Quadratic}, $\Alt(4)$ is abelian, a contradiction.

Thus, $\dim W = 0$. As $W$ is dim-connected, $2V = W = 0$, so $V$ is an elementary abelian $2$-group. Let $K = \left\{1, \alpha_{1},\alpha_{2},\alpha_{3}\right\}$ be the  Klein $4$-group in $\Alt(4)$. Since every $v\in V$ satisfies $v=-v$, $\ad_{\alpha_i}(v)=(1-\alpha_i)v=(1+\alpha_i)v=\tr_{\alpha_i}(v)$, thus $B_{\alpha_i} = C_{\alpha_i}$. By faithfulness, $B_{\alpha_i} \neq 0$ and $C_{\alpha_i} \ne V$, so $B_{\alpha_i} = C_{\alpha_i}$ must have dimension $1$ for each $i$. 

We claim that $B_{\alpha_i} = B_{\alpha_j}$ for all $i,j\in \{1,2,3\}$. All elements of $K$ commute, so $K$ acts on each $B_{\alpha_i}$ (see Lemma~\ref{l.PermutingBandC}). Thus, 
 $\ad_{\alpha_i}(B_{\alpha_j})\le B_{\alpha_j}$. If $\ad_{\alpha_i}(B_{\alpha_j})$ has dimension $1$, then $B_{\alpha_j}=\ad_{\alpha_i}(B_{\alpha_j})\le\ad_{\alpha_i}(V)=B_{\alpha_i}$, so  $B_{\alpha_j}=B_{\alpha_i}$ (as both groups are dim-connected of dimension $1$). 
 If instead $\ad_{\alpha_i}(B_{\alpha_j})$ has dimension $0$, then $\ad_{\alpha_i}(B_{\alpha_j})=0$, so $\alpha_i$ fixes $B_{\alpha_j}$. Thus, $B_{\alpha_j} \le C_V(\alpha_i)$.
Since we  know each of $B_{\alpha_j}$, $B_{\alpha_i}$, and $C_V(\alpha_i)$ have dimension $1$, $B_{\alpha_j}$ and $B_{\alpha_i}$ must both be the dc-component of $C_V(\alpha_i)$, which is unique (see Remark~\ref{rem.dcComponent}). So, in all cases, $B_{\alpha_i} = B_{\alpha_j}$.

Since the $K$ is normal in $\Sym(4)$ and $B_{\alpha_i} = B_{\alpha_1}$ for all $i\in \{1,2,3\}$, Lemma~\ref{l.PermutingBandC} shows that $B_{\alpha_1}$ is a $\Sym(4)$-submodule of $V$. Considering the action of $\Sym(4)$ on $B_{\alpha_1}$ and $V/B_{\alpha_1}$, we find (as before) that $\Alt(4)$ acts quadratically on $V$, implying that $\Alt(4)$ is abelian, a contradiction.
\end{proof}


\begin{lemma}\label{l.minimalAnIrreducible}
Assume $n\ge 5$. Let $V \in \Mod(\Alt(n), d)$ be faithful. If $d$ is the minimal dimension of a faithful $\Alt(n)$-module, then $V$ is dc-irreducible.
\end{lemma}
\begin{proof}
Let $d$ be the minimal dimension of a faithful $\Alt(n)$-module, and choose $V \in \Mod(\Alt(n), d)$.  Assume $V$ is \emph{not} dc-irreducible, and let $W$ be a proper nontrivial dim-connected submodule of $V$. Then $W$ and $V/W$ are $\Alt(n)$-modules of dimension at most $d-1$, so by assumption, $\Alt(n)$ does not act faithfully on $W$ nor $V/W$. Since $n\ge 5$, $\Alt(n)$ is simple, so in fact, $\Alt(n)$ acts trivially  $W$ and $V/W$, meaning that $\Alt(n)$ acts quadratically on $V$. By Lemma~\ref{l.Quadratic}, $\Alt(n)$ is abelian, a contradiction.
\end{proof}

\begin{lemma}\label{l.Alt5}
If $V \in \Mod(\Alt(5), d)$ is faithful, then $d \ge 2$. Further, $d=2$ is only possible when $V$ has characteristic $2$. 
\end{lemma}
\begin{proof}
 Since $\Alt(5)$ contains $\Alt(4)$, the first part of the statement is just Lemma~\ref{l.Alt4}. Now assume $d=2$. Then again by Lemma~\ref{l.Alt4}, $V$ is \emph{not} $2$-divisible, so the submodule $2V$ is proper in $V$. Since $2V$ is dim-connected, Lemma~\ref{l.minimalAnIrreducible} implies that $2V = 0$, so $V$ has characteristic $2$.
\end{proof}

\begin{remark}
    Note that there do exist faithful $V \in \Mod(\Alt(5), 2, 2)$ since $\Alt(5) \cong \SL_2(\bF_4)$, which acts naturally on the $2$-dimensional $\bF_4^2$.
\end{remark}

\begin{lemma}\label{l.Alt6}
If $V \in \Mod(\Alt(6), d)$ is faithful, then $d \ge 3$. Further, $d=3$ is only possible when $V$ has characteristic $3$, and $d=4$ is only possible when $V$ has characteristic $2$ or $3$.
\end{lemma}
\begin{proof}
By Lemma~\ref{l.minimalAnIrreducible} and Fact~\ref{fact.Characteristic}, $V$ has a characteristic $q$, and Fact~\ref{fact.Divisibility} implies that $V$ $p$-divisible for every prime $p$ not equal to $q$.

\begin{claim}\label{l.Alt6:cl:dAtLeast3}
$d \ge 3$.
\end{claim}
\begin{proofclaim}
Consider the subgroup $S$ of all permutations in $\Alt(6)$ that either fix both $5$ and $6$ or interchange $5$ and $6$:
\[
S=\left\{ \rho\in\Alt(6)\mid \{\rho(5),\rho(6)\} = \{5,6\}\right\}.
\]
It can be checked that  $S\cong \Sym(4)$, so by Lemma~\ref{l.Sym4}, $d\ge 3$.
\end{proofclaim}



\begin{claim}\label{l.Alt6:cl:dichotomy}
Suppose $q\ne 3$ and $d \le 4$. If $\gamma,\delta \in \Alt(6)$ are disjoint $3$-cycles, then either 
\begin{itemize}
    \item $B_\gamma = B_\delta$ and $C_\gamma = C_\delta$, or
    \item  $B_\gamma = C_\delta$ and $C_\gamma = B_\delta$.
\end{itemize} 
\end{claim}
\begin{proofclaim}
As $\gamma$ and $\delta$ commmute, $\langle\delta\rangle$ acts on $B_\gamma$ by Lemma~\ref{l.PermutingBandC}\eqref{i:BCConjugation}, so we apply Fact~\ref{fact.Coprimality}  to see $B_\gamma = \ad_\delta(B_\gamma) \qoplus \tr_\delta(B_\gamma)$. Let $W_1 := \ad_\delta(B_\gamma)$ and $W_2 := \tr_\delta(B_\gamma)$. 
 
We first show that neither $W_1$ nor $W_2$ has dimension $1$. We may assume $\gamma = (123)$ and $\delta = (456)$. Noting that $(12)(45)\delta(12)(45) = \delta^{-1}$, we find that  $\Sigma = \langle \delta, (12)(45)\rangle$ is isomorphic to $\Sym(3)$. Every element of $\Sigma$ centralizes or inverts $\gamma$, so Lemma~\ref{l.PermutingBandC} implies that $\Sigma$ acts on $B_{\gamma}$. Thus, we may apply Corollary~\ref{c.Bgamma} to see $\dim W_1 = \dim \ad_\delta(B_{\gamma}) \neq 1$. A similar argument shows that $\dim \ad_\gamma(C_{\delta}) \neq 1$, and as $\gamma$ and $\delta$ commute, $\ad_\gamma(C_{\delta}) = (1-\gamma)(1+\delta+\delta^2)(V) = (1+\delta+\delta^2)(1-\gamma)(V) = \tr_\delta(B_{\gamma}) = W_2$.
 
We next show $W_1$ and $W_2$ cannot both have dimension larger than $1$. Suppose they do. 
In this case, $\dim B_\gamma \ge 4$. However, $W_2=\tr_\delta(B_\gamma) \le  C_\delta$, so we also find that $\dim C_\delta \ge 2$.
As $\gamma$ and $\delta$ are conjugate, $\dim C_\gamma \ge 2$, so Fact~\ref{fact.Coprimality} implies that $\dim V =  \dim B_\gamma + \dim C_\gamma \ge 4+2$, a contradiction. 
Thus $\dim W_i = 0$ for some $i$. 

Suppose $\dim W_1 = 0$; then $B_\gamma = W_2 \le \tr_\delta(V) =  C_\delta$. As the action is faithful, Lemma~\ref{l.AdjointTrace}\eqref{l.AdjointTrace.i:BTrivial} implies that $W_2 = B_\gamma \neq 0$, so as $\dim W_2 \neq 1$, it must be that $\dim W_2 \ge 2$. Thus,  $\dim C_\delta\ge 2$, and as $\gamma$ and $\delta$ are conjugate, $\dim C_\gamma \ge 2$. Since
$\dim B_\gamma + \dim C_\gamma=\dim V \le 4$, it must be that $\dim B_\gamma = \dim C_\gamma = \dim C_\delta= 2$. By dim-connectedness of $C_\delta$, $B_\gamma = C_\delta$, and by conjugacy, $B_\delta = C_\gamma$.

Next suppose $\dim W_2 = 0$; then $B_\gamma = W_1 \le \ad_\delta(V) = B_\delta$, so $B_\gamma \le B_\delta$. Since $\gamma$ and $\delta$ are conjugate, we find $B_\gamma = B_\delta$. 
Now, $\ad_\delta(C_\gamma) \le B_\delta = B_\gamma$, and as $\delta$ acts on $C_\gamma$ (since $\delta$ commutes with $\gamma$), we also have that $\ad_\delta(C_\gamma) \le C_\gamma$. By Fact~\ref{fact.Coprimality}, $\ad_\delta(C_\gamma)$, which is contained in  $B_\gamma\cap C_\gamma$, has dimension $0$. Another application of Fact~\ref{fact.Coprimality} shows $C_\gamma = \ad_\delta(C_\gamma)\qoplus \tr_\delta(C_\gamma)$, so $\dim C_\gamma = \dim \tr_\delta(C_\gamma)$. By dim-connectedness of $C_\gamma$, $C_\gamma = \tr_\delta(C_\gamma)$, so  $C_\gamma\le C_\delta$. By conjugacy, $C_\gamma = C_\delta$.
\end{proofclaim}

\begin{claim}\label{l.Alt6:cl:outer}
Suppose $q\ne 3$ and $d \le 4$.  Up to composing the action  on $V$ with an automorphism of $\Alt(6)$, $B_\gamma = B_\delta$ for any  disjoint $3$-cycles $\gamma,\delta \in \Alt(6)$. 
\end{claim}
\begin{proofclaim}
By Claim~\ref{l.Alt6:cl:dichotomy}, we only need to consider when $B_\gamma = C_\delta$ and $C_\gamma = B_\delta$.
Let $\varepsilon = \gamma \delta$. We prove $B_\varepsilon = V$.
As $B_\gamma = C_\delta\le C_V(\delta)$, $\ad_\varepsilon(B_\gamma) = (1-\gamma\delta)(B_\gamma) = (1-\gamma)(B_\gamma) = \ad_\gamma(B_\gamma)$. Since $V$ is $3$-divisible, Lemma~\ref{l.AdjointTrace}\eqref{l.AdjointTrace.i:BCIntersection} implies that  $B_\gamma \cap C_V(\gamma)$ has dimension $0$. Now, $B_\gamma \cap C_V(\gamma)$ is the kernel of the map $\ad_\gamma : B_\gamma \rightarrow B_\gamma$, so the image of the map has dimension equal to $\dim B_\gamma$. Thus, $B_\gamma = \ad_\gamma(B_\gamma) = \ad_\varepsilon(B_\gamma)\le B_\varepsilon$. A similar argument shows  $B_\delta  \le B_\varepsilon$. We are assuming that $C_\gamma = B_\delta$, so $B_\gamma + C_\gamma \le B_\varepsilon$. By Fact~\ref{fact.Coprimality}, $V = B_\gamma + C_\gamma$, so $B_\varepsilon = V$ as desired.

This also holds of any conjugate $\varepsilon'$ of $\varepsilon$.
Now, there is an (outer) automorphism of $\Alt(6)$ that maps the conjugacy class of $\gamma$  to the conjugacy class of $\varepsilon$, see for example \cite{JR82}.
Up to composing with this automorphism, we now have $B_\gamma = V$. But this also holds of any conjugate of $\gamma$, and $\delta$ is one such. So we returned to $B_\gamma = B_\delta$ and we are done.
\end{proofclaim}

\begin{claim}\label{l.Alt6:cl:dAtLeast4}
If $q\ne 3$, then $d \ge 4$.  
\end{claim}
\begin{proofclaim}
By Claim~\ref{l.Alt6:cl:dAtLeast3}, $d \ge 3$, so assume $d = 3$. Let $\gamma,\delta \in \Alt(6)$ be disjoint $3$-cycles. Using Claim~\ref{l.Alt6:cl:outer}, we consider $V$ as an $\Alt(6)$-module for which $B_\gamma = B_\delta$. Set $B = B_\gamma = B_\delta$ and  $\varepsilon = \gamma\delta$.

We first show $\ad_\delta(B) = B$. Indeed, $\tr_\delta(B) = \tr_\delta(B_\delta) \leq B_\delta \cap C_V(\delta)$ has characteristic $3$ by Lemma~\ref{l.AdjointTrace}\eqref{l.AdjointTrace.i:BCIntersection}, so it has dimension $0$. We conclude by coprimality.

We then show that $B_\varepsilon = B$. Note that:
\begin{align*}
    B_\varepsilon 
    = (1-\gamma\delta)(V)
    &= [(1-\delta) + (1-\gamma)\delta](V)\\
    &\le (1-\delta)(V) + (1-\gamma)\delta(V)
    \le B_\delta + B_\gamma = B + B = B.
\end{align*} 
Now consider the decomposition $B = \ad_\varepsilon(B) \qoplus \tr_\varepsilon(B)$. As in the proof of Claim~\ref{l.Alt6:cl:dichotomy}, we find that  $\ad_\varepsilon(B)$ and $\tr_\varepsilon(B)$ are both invariant under the action of $\langle \gamma, (12)(45)\rangle \cong \Sym(3)$, so by Corollary~\ref{c.Bgamma}, neither have dimension $1$. Since $\dim B \le 3$, one of $\ad_\varepsilon(B)$ or $\tr_\varepsilon(B)$ must have dimension $0$. If $\dim \ad_\varepsilon(B) = 0$, then $B = \tr_\varepsilon(B) \le C_\varepsilon$, so $B_\varepsilon \le B \le C_\varepsilon \le C_V(\varepsilon)$. By Lemma~\ref{l.AdjointTrace}\eqref{l.AdjointTrace.i:BCIntersection}, we find that $B_\varepsilon$ has characteristic $3$, which contradicts the fact that $V$ is $3$-divisible. Thus $\dim \tr_\varepsilon(B) = 0$, so  
 $B = \ad_\varepsilon(B) \le B_\varepsilon\le B$, implying  $B= B_\varepsilon$.
 
We  next aim  to show that $\gamma$ and $\delta$ agree on $B$. By Lemma~\ref{l.AdjointTrace}\eqref{l.AdjointTrace.i:traceonB}, $\tr_\gamma(B) = \tr_\delta(B) = \tr_\varepsilon(B) = 0$, so, \emph{restricted to $B$}, we have the following equations.
\[\begin{array}{ccccccc}
    1 & + & \gamma & + & \gamma^2 & = & 0\\
    1 & + & \delta & + & \delta^2 & = & 0\\
    1 & + & \gamma\delta & + & \gamma^2\delta^2 & = & 0
\end{array}\]
Multiplying the third equation by $\gamma$ and subtracting it from the first yields
\[
    0 = 1 + \gamma^2 - \gamma^2\delta - \delta^2
    = (1-\delta^2) + \gamma^2(1-\delta)
    = (1+\delta + \gamma^2)(1-\delta).
\]
Then, using that $1 + \delta = - \delta^2$ (from the second equation), we obtain 
\[(\gamma^2 -\delta^2)(1-\delta) = 0.\]
This equation holds on $B$, so we have $(\gamma^2 -\delta^2)\circ\ad_\delta(B) = 0$.
But $\ad_\delta (B) = B$,
so $\gamma^2 = \delta^2$ on $B$. Squaring both sides, we find that $\gamma=\delta$ on $B$.

We also have  $C_\gamma = C_\delta$ by Claim~\ref{l.Alt6:cl:dichotomy}.  Setting $C= C_\gamma = C_\delta$, we see that $\gamma=\delta$ on $C$ since both maps centralize $C$ by Lemma~\ref{l.AdjointTrace}\eqref{l.AdjointTrace.i:CCentral}. Since $V=B + C$, we now have that $\gamma=\delta$ on $V$, which is a contradiction to faithfulness.
\end{proofclaim}

\begin{claim}\label{l.Alt6:cl:dAtLeast5}
If $q\ne 2$ and $q\ne 3$, then $d \ge 5$.  
\end{claim}
\begin{proofclaim}
From our previous work, it remains to consider when $d = 4$. Since $V$ is $2$-divisible, Lemma~\ref{l.WeightLemma} shows that $\dim V=3\cdot\dim V_{\alpha}+\dim C_{K}$ where $K\le \Alt(4)$ is the Klein $4$-group on $\{1, 2, 3, 4\}$ and $\alpha$ is any nontrivial element of $K$. As $\dim V = 4$, $\dim C_{K}$ is either $1$ or $4$. If $\dim C_{K} = 4$, then $V = C_{K} \le C_V(K)$, contradicting faithfulness, so  $\dim C_{K} = 1$. 

Now, the permutations $(13)(56)$ and $(12)(56)$ normalize $K$, so they act on $C_{K}$.
Notice that they are conjugate under $(123)$, which normalizes $K$ as well.
By Remark~\ref{r.dimension1}, the product $(123) = (13)(56)(12)(56)$ centralizes $C_{K}$.
Then $\ad_{(123)}(C_{K}) = 0$, so Fact~\ref{fact.Coprimality} implies that $\tr_{(123)}(C_{K}) = C_{K}$. A similar argument shows that $\tr_{\gamma}(C_{K}) = C_{K}$ for any $3$-cycle in $\Alt(4)$, so in particular, $C_K \le C_{(123)} \cap C_{(234)}$. 

By Claims~\ref{l.Alt6:cl:dichotomy} and \ref{l.Alt6:cl:outer}, we may assume  $C_{(123)} = C_{(456)}$ and $C_{(234)} = C_{(156)}$. Thus, $C_K \le C_{(123)} \cap C_{(234)} \cap C_{(456)} \cap C_{(156)}$, so by Lemma~\ref{l.AdjointTrace}\eqref{l.AdjointTrace.i:CCentral}, $C_K$ is centralized by each of $(123),(234),(456),(156)$. We conclude that $C_K \le C_V(\langle (123),(234),(456),(156) \rangle)$. Now, direct computations show  $ \Alt(6)= \langle (123),(234),(456),(156) \rangle $, so $C_K$ is centralized by all of $\Alt(6)$. In particular, $C_K$ is a proper nontrivial $\Alt(6)$-submodule of $V$, so $V$ is \emph{not} dc-irreducible. However, our assumptions together with the result of Claim~\ref{l.Alt6:cl:dAtLeast4} imply that  $d=4$ is the minimal dimension of a faithful $\Alt(6)$-module, so we have a contradiction to Lemma~\ref{l.minimalAnIrreducible}.
\end{proofclaim}
This completes the proof of Lemma~\ref{l.Alt6}.
\end{proof}\setcounter{claim}{0}

We have everything needed for Theorem~\ref{t.AltBounds}.

\begin{proof}[Proof of Theorem~\ref{t.AltBounds}]
Lemmas~\ref{l.Alt4}, \ref{l.Alt5},  \ref{l.Alt6} prove almost everything. The only point needing clarification is when $V\in\Mod(\Alt(4),2)$. In this case, Lemma~\ref{l.Alt4} implies that $V$ is \emph{not} $2$-divisible, but this is equivalent to $\Omega_2(V) = \{v \in V: 2v = 0\}$ having dimension at least $1$ by Fact~\ref{fact.Divisibility}.
\end{proof}

\begin{remark}\label{r.ProofOfTable}
We briefly discuss how the right side of Table~\ref{tab.SmallDimension} now follows. 
Letting $d$ be minimal such that $\Mod(\Alt(n),d,p)$ contains a faithful module, Theorem~\ref{t.AltBounds} and the main result of \cite{CDW23} combine to show that $d$ is no smaller than what is stated in Table~\ref{tab.SmallDimension}. (The exceptional cases for $n=7$ and $n=8$ are implied by the lower bound when $n=6$.) That $d$ is no larger than stated follows from the examples in \S\S~\ref{ss.FamiliarModules}, except when $n=7,8$. The upper bound on $d$ in these remaining cases comes from the natural action of $\Alt(8)\cong \SL(4,2)$ on a $4$-dimensional vector space of characteristic $2$ (discussed in \S\S~\ref{ss.RemainingWork} following Task~\ref{p.AltInChar2}).
\end{remark}

\section{Recognizing minimal modules}

We now work to identify various minimal modules; this will culminate in  a proof of Theorem~\ref{t.AltRec}.

\subsection{Minimal modules for \texorpdfstring{$\Sym(3)$}{Sym(3)} and \texorpdfstring{$\Sym(4)$}{Sym(4)}}
We begin with recognition results for the minimal, faithful $\Sym(n)$-modules in the exceptionally small cases of $n=3,4$. The main tool is the Recognition Lemma of \cite{CDW23}, and to apply it, we mainly just need to show that a minimal, faithful module is also dc-irreducible (with appropriate restriction on the characteristic). However, our work here exposed a small issue with the Recognition Lemma in the case when $n=4$, which we first discuss and correct.

\begin{remark}
The Recognition Lemma of \cite{CDW23} is incorrect for $n=4$. To see why, 
first note that $\std(4,\mathbb{C}^+)$ and $\std^\sigma(4,\mathbb{C}^+)$ are both faithful and irreducible $\Sym(4)$-modules of dimension $3$; let $V$ denote either one of these modules.  Now, the condition ``$\Alt(\{3,4\})$ centralizes $\ad_{(12)}(V)$'' is trivially satisfied (for either choice of $V$) because $\Alt(\{3,4\})$ is itself trivial; this is the issue with $n=4$. Consequently, the Recognition Lemma implies that $V$ must be isomorphic to $\std(4,\mathbb{C}^+)$, but $\std(4,\mathbb{C}^+)\not\cong \std^\sigma(4,\mathbb{C}^+)$. 

When $n=3$, the condition ``$\Alt(\{3\})$ centralizes $\ad_{(12)}(V)$'' also holds trivially, but here we are saved by the fact that $\std(3,\mathbb{C}^+)\cong\std^\sigma(3,\mathbb{C}^+)$.
\end{remark}

We now provide a corrected version of the Recognition Lemma, the proof of which assumes familiarity with (and directly references) the original proof.

\begin{lemma}[Corrected version of {\cite[Recognition Lemma]{CDW23}}]\label{lemma.RecognitionUpdate}
If $n\neq 4$, the Recognition Lemma of \cite{CDW23} holds as stated; if $n= 4$, the Recognition Lemma holds under the following additional assumption: for any transposition $\tau$, $B_\tau$ is either centralized or inverted by each element of $\Sym(\tau^\perp)$, where $\tau^\perp$ is the set of fixed points of $\tau$. 

In particular, the Recognition Lemma holds when  $\dim B_\tau = 1$ for $\tau$ any transposition $\tau$.
\end{lemma}
\begin{proof}
We adopt the notation of \cite{CDW23}; set $S = \Sym(n)$ and $A = S' = \Alt(n)$.
Studying the proof of \cite[Recognition Lemma]{CDW23}, there is only one error; it occurs in the proof of Claim~1 when showing that $B_{(ij)} \cap B_{(jk)}$ is $S$-invariant. We proceed to review the proof that $B_{(ij)} \cap B_{(jk)}$ is $S$-invariant and correct the issue, which only appears when $n=4$.

In \cite{CDW23}, it is first claimed that $B_{(ij)} \cap B_{(jk)}$ is $S_{(ijk)}$-invariant, where $S_{(ijk)} = \Sym(\{i,j,k\})$. There is no error with this, and in particular, this shows that $B_{(ij)} \cap B_{(jk)}$ is $S$-invariant in the case when $n=3$. 

The authors next claim that $B_{(ij)} \cap B_{(jk)}$ is $A_{j^\perp}$-invariant, where $A_{j^\perp}$ denotes the subgroup of $A$ consisting of those permutations fixing $j$. Their reasoning is that $A_{(ij)^\perp}$ (defined to be $\Alt((ij)^\perp)$) centralizes $B_{(ij)}$ and $A_{(jk)^\perp}$  centralizes $B_{(jk)}$, so $A_{j^\perp}$ centralizes $B_{(ij)} \cap B_{(jk)}$ since $A_{j^\perp} = \langle A_{(ij)^\perp}, A_{(jk)^\perp} \rangle$. Now, $A_{(ij)^\perp}$ does indeed centralizes $B_{(ij)}$ (and similarly for $A_{(jk)^\perp}$) by the assumptions  of the Recognition Lemma, so $B_{(ij)} \cap B_{(jk)}$ is centralized by $\langle A_{(ij)^\perp}, A_{(jk)^\perp} \rangle$. However---and here is the error---we only have that $A_{j^\perp} = \langle A_{(ij)^\perp}, A_{(jk)^\perp} \rangle$ when $n\ge 5$. If $n\ge 5$, we may fix $a,b\notin \{ijk\}$, and note that $\langle A_{(ij)^\perp}, A_{(jk)^\perp} \rangle$ contains all three cycles of the form $(abc)$ for fixed $a,b\notin \{ijk\}$ and variable $c\neq j$, which do indeed generate $A_{j^\perp}$. But if $n=4$, $A_{(ij)^\perp}=A_{(jk)^\perp}=1$ while $A_{j^\perp}$ is nontrivial.

Now, assuming $n=4$,  we show $B_{(ij)} \cap B_{(jk)}$ is $A_{j^\perp}$-invariant under the (stronger) assumptions of the present lemma. Here we have that $B_{(ij)}$ is either centralized or inverted by each element of $S_{(ij)^\perp}$ and similarly for $(jk)$. Thus, \emph{every} subgroup of $B_{(ij)}$ is $S_{(ij)^\perp}$-invariant, and  \emph{every} subgroup of $B_{(jk)}$ is $S_{(jk)^\perp}$-invariant. We conclude that $B_{(ij)} \cap B_{(jk)}$ is invariant under the action of $\langle S_{(ij)^\perp}, S_{(jk)^\perp} \rangle$, and, crucially,  it is easy to verify that $S_{j^\perp} = \langle S_{(ij)^\perp}, S_{(jk)^\perp} \rangle$. Thus,  $B_{(ij)} \cap B_{(jk)}$ is invariant under $S_{j^\perp}$, hence also $A_{j^\perp}$.

Finally, since $B_{(ij)} \cap B_{(jk)}$ is both $S_{(ijk)}$- and $A_{j^\perp}$-invariant, we find (as in \cite{CDW23}) that $B_{(ij)} \cap B_{(jk)}$ is $S$-invariant since $S = \langle S_{(ijk)}, A_{j^\perp} \rangle$.

We now address the special case when $\dim B_\tau = 1$ for $\tau$ a transposition. Under this assumption, Lemma~\ref{l.Sym2} shows that $B_\tau$ is either centralized or inverted by a transposition of $S_{\tau^\perp}$, so the same holds of all element in $S_{\tau^\perp}$. Since  transpositions in $S_{\tau^\perp}$ are conjugate, they all act the same on $B_\tau$, so $B_\tau$ is centralized by $A_{\tau^\perp}$ (as required by the original  Recognition Lemma).
\end{proof}

\begin{remark}
The remainder of \cite{CDW23} needs no further correction as the Recognition Lemma is only used when $n\ge 7$; this is fairly easy to track as \cite{CDW23} is structured so that the Recognition Lemma is essentially only used in the proof of the main theorem in the final section (though it also appears in a remark immediately preceding the proof of the Theorem).
\end{remark}

We have yet to give a proper statement of the Recognition Lemma from \cite{CDW23}, which we will need. This comes next, but we first introduce  terminology (adapted from the algebraic context) to streamline the presentation.

\begin{definition}
Let $V$ and $W$ be modules in some modular universe with an additive dimension. A compatible homomorphism $\varphi:V\rightarrow W$ is called a \textbf{dim-isogeny} (for ``dimensional isogeny'') if it is surjective with  kernel of dimension $0$. 
\end{definition}

We now give a  simplified version of \cite[Recognition Lemma]{CDW23}, which incorporates the correction of Lemma~\ref{lemma.RecognitionUpdate}.

\begin{fact}[Special case of {\cite[Recognition Lemma]{CDW23}}]\label{fact.Recognition}
Let $n\ge 2$. Suppose that $V \in \Mod (\Sym(n), d, p)$ is faithful and dc-irreducible; assume $p$ does not divide $n$. 
If $\dim B_{(12)} = 1$, then there is a $\Sym(n)$-module dim-isogeny  $\varphi:\std(n,L) \rightarrow V$ for some $1$-dimensional subgroup $L\le V$. Further, if $V$ has no elements of order dividing $n$, then $\varphi$ is an isomorphism.
\end{fact}

\begin{lemma}\label{l.Sym3Recognition}
Let $V\in \Mod(\Sym(3),2)$ be faithful and $3$-divisible. Then, there is a $\Sym(3)$-module dim-isogeny $\varphi:\std(3,L)\rightarrow V$ for some $1$-dimensional subgroup $L\le V$; further, if $V$ has no elements of order $3$, then $\varphi$ is an isomorphism.
\end{lemma}
\begin{proof}
Assume $V$ is as stated in the lemma; we aim to apply Fact~\ref{fact.Recognition}.  

To see that $V$ is dc-irreducible, assume,  towards a contradiction, that $W$ is a proper nontrivial dim-connected $\Sym(3)$-submodule of $V$. Then $\dim W=1$.  
Thus, $W$ and $V/W$ are both $1$-dimensional $\Sym(3)$-modules. By Lemma~\ref{l.Sym3}, $\Sym(3)$ does \emph{not} act faithfully on $V/W$, so the $3$-cycles act trivially on $V/W$. Let $\gamma\in \Sym(3)$ be a $3$-cycle. Since $\gamma$ centralizes $V/W$, 
one has $\ad_\gamma (V/W) = 0$ so $\ad_\gamma(V) \le W$.
By faithfulness of $\Sym(3)$ on $V$, $\ad_\gamma(V) \neq 0$, so as $W$ is dim-connected of dimension $1$, it must be that $\ad_\gamma(V) = W$, which contradicts Corollary~\ref{c.Bgamma}.

To apply Fact~\ref{fact.Recognition}, it remains to show  $\dim B_{(12)} = 1$. Of course, $\dim B_{(12)} \le \dim V = 2$. If  $\dim B_{(12)} = 0$, then $(12)$ centralizes $V$, contradicting faithfulness. So suppose $\dim B_{(12)} = 2$. Then $B_{(12)} = V$, so $(12)$ inverts $V$. By Lemma~\ref{l.PermutingBandC}\eqref{i:BCConjugation}, $(23)$ inverts $V$, so $(123) = (12)(23)$ centralizes $V$, again contradicting faithfulness. We conclude that $\dim B_{(12)} = 1$.
\end{proof}

\begin{lemma}\label{l.Sym4Recognition}
Let $V\in \Mod(\Sym(4),3)$ be faithful and $2$-divisible. Then, there is a $\Sym(4)$-module dim-isogeny  $\varphi:\hat V\rightarrow V$ with domain $\std(4,L)$ or $\std^\sigma(4,L)$ for some $1$-dimensional subgroup $L\le V$;  further, if $V$ has no elements of order $2$, then $\varphi$ is an isomorphism.
\end{lemma}
\begin{proof}
Assume $V$ is as stated in the lemma; we again use Fact~\ref{fact.Recognition}.

We first  show  $V$ is dc-irreducible. If $W$ is a proper nontrivial dim-connected $\Sym(4)$-submodule of $V$, then $W$ and $V/W$ are dim-connected, $\Sym(4)$-modules of dimension at most $2$. By Lemma~\ref{l.Sym4}, $\Sym(4)$ does not act faithfully on $W$ nor $V/W$. Thus, the Klein $4$-group $K$ acts trivially on both $W$ and $V/W$. If we choose any nontrivial $\alpha \in K$, then (as in the proof of Lemma~\ref{l.Sym3Recognition}) $\ad_\alpha(V) \le W$, so   $\ad_\alpha(\ad_\alpha(V)) \le\ad_\alpha(W) =0$. Thus, $\ad_\alpha(V) \le C_V(\alpha)$ by Lemma~\ref{l.AdjointTrace}\eqref{l.AdjointTrace.i:BTrivial}, and then Lemma~\ref{l.AdjointTrace}\eqref{l.AdjointTrace.i:BCIntersection} implies that $\ad_\alpha(V)$ is an elementary abelian $2$-group. Since $V$ is $2$-divisible, $\ad_\alpha(V)=0$ by Fact~\ref{fact.Divisibility}, so $\alpha$ acts trivially on $V$, a contradiction. We conclude that $V$ is dc-irreducible.

Next, exactly as in  Lemma~\ref{l.Sym3Recognition}, we find that $B_{(12)}$ does not have dimension $0$ nor $3$, and if $\dim B_{(12)} = 1$, we directly apply Fact~\ref{fact.Recognition} to get the desired result with $\hat V=\std(4,L)$. 

Now assume $\dim B_{(12)} = 2$; by Fact~\ref{fact.Coprimality},  $\dim C_{(12)} = 1$. In this case, we aim for the conclusion of the lemma in which $\hat V=\std^\sigma(4,L)$.
Consider the module $V^\sigma$ with the same underlying group $V$ and $\Sym(4)$-action defined by $g*v = \sgn(g)(g\cdot v)$, where $g\cdot v$ is the original action of $g$ on $v$. Then, $\ad_{(12)}(V^\sigma) = \tr_{(12)}(V)$, so $\dim \ad_{(12)}(V^\sigma) = \dim C_{(12)} = 1$. Thus, Fact~\ref{fact.Recognition} applies to $V^\sigma$, and there is  a surjective $\Sym(4)$-module homomorphism $\varphi^\sigma:\std(4,L)\rightarrow V^\sigma$ with the required properties. This then induces a $\Sym(4)$-homomorphism $\varphi:\sgn(4,\bZ)\otimes\std(4,L)\rightarrow \sgn(4,\bZ)\otimes V^\sigma$, so as $\sgn(4,\bZ)\otimes V^\sigma \cong V$, we have our desired morphism (found by composing $\varphi$ with the appropriate isomorphism).
\end{proof}

\subsection{Minimal modules for \texorpdfstring{$\Alt(5)$}{Alt(5)} in characteristic 2}

Here we study $V\in \Mod(\Alt(5),2,2)$. In this setting, we have $\Alt(5)\cong \SL(2,\bF_4)$, and we aim to identify $V$ as the natural $\SL(2,\bF_4)$-module (see Example~\ref{eg.NatSL2}). The main tool is the following (lovely) result of Timmesfeld; see also \cite{SmS89}.

\begin{fact}[Special case of {\cite[(3.7)~Corollary]{TiF01}}]\label{fact.Nat(SL2)}
Assume $G\cong \SL(2,\bF)$ or $G\cong \PSL(2,\bF)$ for $\bF$ a field, and let $U\le G$ be the image of the subgroup of strictly upper triangular matrices. Suppose  $V$ is a faithful $G$-module such that
\begin{enumerate}
    \item $U$ acts quadratically on $V$;
    \item $[G,V] = V$ and $C_V(G) = 0$.
\end{enumerate}
Then, $G\cong \SL(2,\bF)$, and for some  $I$, $V= \bigoplus_{i\in I} V_i$ with each $V_i \cong \Nat(2,\bF)$.
\end{fact}

\begin{lemma}\label{l.Alt5Recognition}
    Let $V\in \Mod(\Alt(5),2,2)$ be faithful. Then $V$ is an $\bF_4$-vector space, and viewing $\Alt(5)$ as $\SL(2,\bF_4)$, we have $V\cong \Nat(2, L)$ for some $1$-dimensional subgroup $L\le V$.
\end{lemma}
\begin{proof}
Let $V$ be as stated, and view $V$ as a  $G$-module for $G =  \SL(2,\bF_4)$.
By Lemmas~\ref{l.minimalAnIrreducible} and \ref{l.Alt5}, $V$ is dc-irreducible. 

\begin{claim}
For some index set $I$, $V= \bigoplus_{i\in I} V_i$ with each $V_i \cong \Nat(2,\bF_4)$.
\end{claim}
\begin{proofclaim}
To apply Fact~\ref{fact.Nat(SL2)}, we first aim to show that the subgroup $U$ of strictly upper triangular matrices  acts quadratically on $V$. Note that $U$ is an elementary abelian $2$-group of order $4$. By Lemma~\ref{l.primality}, $[U,V] < V$, and since the action is faithful, $[U,V] \neq 0$. Additionally, $V$ is dim-connected of dimension $2$, so $[U,V] < V$ implies $\dim [U,V] \le 1$. Applying Lemma~\ref{l.primality} again, we have that $[U,U,V] < [U,V]$. Since $[U,V]$ is also dim-connected,  $\dim [U,U,V] =0$, and as  $[U,U,V]$ is dim-connected, $[U,U,V] = 0$. 

We next show that  $[G,V] = V$ and $C_V(G) = 0$. The first follows immediately from the fact that $V$ is dc-irreducible, so it remains to show $C_V(G)=0$.

As in Claim~\ref{l.Alt6:cl:dAtLeast3} of Lemma~\ref{l.Alt6}, $\Alt(5)$ (hence $G$) contains a subgroup $S$ isomorphic to $\Sym(3)$. Viewing $V$ as an $S$-module, Lemma~\ref{l.Sym3Recognition} tells us that $V$ is isomorphic $\std(3,L)$ for some $1$-dimensional group $L$. In the current setting, $\std(3,L)$, hence $L$, has characteristic $2$, and as $2$ does not divide $3$, no element of $\std(3,L)$ is fixed by all of $S$ (also see the Remarks on page 4 of \cite{CDW23}). Thus  $C_V(S)$, hence $C_V(G)$, is trivial.
\end{proofclaim}

\begin{claim}\label{l.Alt5Recognition:cl:Nat(2,L)}
$V\cong \Nat(2, L)$ for $L = [U,V]= C_V(U)$.
\end{claim}
\begin{proofclaim}
We now have $V= \bigoplus_{i\in I} V_i$ with each $V_i \cong \Nat(2,\bF_4)$.
Set $L = [U,V]$. 
It is easily verified that $L = C_V(U)$ and,
more precisely, that $L = \bigoplus_I C_{V_i}(U) = \bigoplus_I L_i$ where $L_i = C_{V_i}(U) \cong(\bF_4)_+$.
We equip $L$ and each $L_i$ with the $G$-trivial action; then $V_i \cong \Nat(2, \bF_4) \otimes L_i$. Since direct sums and tensor products commute, we then have:
\[V = \bigoplus_I V_i \cong \bigoplus_I \Nat(2, \bF_4) \otimes L_i \cong \Nat(2, \bF_4) \otimes \left(\bigoplus_I L_i\right) = \Nat(2, \bF_4) \otimes L,\]
which is the definition of $\Nat(2, L)$.
\end{proofclaim}
This completes the proof.
\end{proof}\setcounter{claim}{0}

\begin{remark}
The equality $C_V(G) = 0$ could also have been obtained as follows: let $\overline{V} = V/C_V(G)$, which has $C_{\overline{V}}(G) = 0$ by perfectness. Proceed to identifying $\overline{V}$. Then basic cohomology proves $C_V(G) = 0$.
\end{remark}

\subsection{Minimal modules for \texorpdfstring{$\Alt(5)$}{Alt(5)} in characteristic 5}

We now investigate $V\in \Mod(\Alt(5),3,5)$. Here we use $\Alt(5)\cong \PSL(2,\bF_5)$ and work to identify $V$ as the adjoint module (as defined in Example~\ref{eg.AdPSL2}). Identification takes place via \emph{cubic} actions, which were first analyzed in the style of Timmesfeld by the second author \cite{DeA16} and later studied quite thoroughly by Gr\"uninger \cite{GrM22}.

\begin{definition}
Let $U$ be a group acting on an abelian group $V$. We say that the action is \textbf{cubic} if $[U,U,U,V] = 0$ but $[U,U,V] \neq 0$.
\end{definition}

\begin{fact}[Special Case of {\cite[Theorem~8.1]{GrM22}}]\label{fact.Ad(PSL2)}
Assume $G\cong \PSL(2,\bF)$ for $\bF$ a field of characteristic not $2$, and let $U\le G$ be the image of the subgroup of strictly upper triangular matrices. Suppose that   $V$ is a faithful $G$-module such that
\begin{enumerate}
    \item $U$ acts cubically on $V$, 
    \item\label{fact.Ad(PSL2).i:C_V(u)} $C_V(u) = C_V(U)$ for all nontrivial $u\in U$,
    \item $[G,V] = V$ and $C_V(G) = 0$, and
    \item $V$ is \emph{square-free}: there does not exist a non-trivial $G$-submodule $W\le V$ for which $U$ acts quadratically on $W$ or on $V/W$.
\end{enumerate}
Then, for some  $I$, $V= \bigoplus_{i\in I} V_i$ with each  $V_i \cong \Ad(\PSL(2,\bF))$.
\end{fact}

The actual statement of \cite[Theorem~8.1]{GrM22} does not assume that the acting group is known to be  $\PSL(2,\bF)$; since we do, we may  remove the square-free condition, as we now show.

\begin{corollary}\label{cor.Ad(PSL2)}
Fact~\ref{fact.Ad(PSL2)} holds without assuming that $V$ is square-free.
\end{corollary}
\begin{proof}
Adopt all hypotheses of Fact~\ref{fact.Ad(PSL2)} with the exception that $V$ is square-free. In particular, we may take $G= \PSL(2,\bF)$.
Suppose $V$ is not square-free, and let $W$ denote the relevant submodule.

First assume that $U$ acts quadratically on $W$. Set $\check W = [G,W]$. As a consequence of $G$ being perfect, $[G,\check W] = \check W$, and we also have $C_{\check W}(G) = 0$ since $C_V(G) = 0$. But then Fact~\ref{fact.Nat(SL2)} applies to the action of $G$ on $\check W$ implying  $G\cong \SL(2,\bF)$ and  contradicting that $G\cong \PSL(2,\bF)$ in characteristic not $2$.

Next, consider when $U$ acts quadratically on $V/W$. 
Let $W_1$ be the preimage in $V$ of $C_{V/W}(G)$, and set $\overline{V} = V/W_1$. We are assuming $[G,V] = V$, so $[G,\overline{V}] = \overline{V}$ as well. We now claim that $C_{\overline{V}}(G) = 0$. Let $W_2$ be the preimage in $V$ of $C_{\overline{V}}(G)$. Then $W \le W_1 \le W_2$ with $[G,W_1]\le W$ and $[G,W_2]\le W_1$. Thus, $[G,G,W_2] \le W$, so as $G$ is perfect, $[G,W_2] \le W$, which implies that $W_2 \le W_1$ and hence $W_2 = W_1$. We may now apply Fact~\ref{fact.Nat(SL2)} to the action of $G$ on $\overline{V}$ and obtain the same contradiction as before. 
\end{proof}

We are looking to identify the faithful modules in $\Mod(\Alt(5),3,5)$ (and later in $\Mod(\Alt(6),3,3)$), which by our previous work are dc-irreducible. Thus, the condition  $[G,V] = V$ is automatically met, but it is not immediate that $C_V(G) = 0$ since dc-irreducibility only forces $\dim C_V(G) = 0$ (and $C_V(G)$ might not be dim-connected). Nevertheless, it is easy to show $C_V(G) = 0$ in the specific cases we care about.
A similar outcome was obtained implicitly in \cite{CDW23}.

\begin{lemma}\label{l.C_V(G)Trivial}
Suppose $\Alt(4) \le G$, and let $V\in \Mod(G,3)$ be faithful.  If $V$ has no elements of order $2$, then $C_V(G) = 0$.
\end{lemma}
\begin{proof}
Let $K=\{1,\alpha_1, \alpha_2, \alpha_3\} \le \Alt(4) \le G$ be a Klein $4$-subgroup. With notation as in Lemma~\ref{l.WeightLemma}, we may write $V = V_{\alpha_1} \qoplus V_{\alpha_2} \qoplus V_{\alpha_3} \qoplus C_K$, and as the $\alpha_i$ are conjugate in $G$, each $V_{\alpha_i}$ has the same dimension. Since $\dim V = 3$, we thus have $V = V_{\alpha_1} \qoplus V_{\alpha_2} \qoplus V_{\alpha_3}$. Now consider $w \in C_V(G)$. Writing $w = v_1+v_2+v_3$ for $v_i\in V_{\alpha_i}$, we may apply  $\alpha_1$ to see that $v_1+v_2+v_3 = \alpha_1(v_1+v_2+v_3) = v_1 - v_2 - v_3$, which implies that $2(v_2+v_3) = 0$. Since $V$ has no elements of order $2$, $v_2+v_3 = 0$. Similarly, using $\alpha_2$ and $\alpha_3$, we find that $v_1+v_3 = v_1+v_2 = 0$. This implies that $2v_i=0$ for each $i$, which, as before, leads to  $v_i = 0$, as desired. 
\end{proof}

We proceed with identification of $V\in \Mod(\Alt(5),3,5)$.

\begin{lemma}\label{l.Alt5RecognitionDim3Char5}
    Let $V\in \Mod(\Alt(5),3,5)$ be faithful. Then $V$ is an $\bF_5$-vector space, and viewing $\Alt(5)$ as $\PSL(2,\bF_5)$, we have $V\cong \Ad(2, L)$ for some $1$-dimensional subgroup $L\le V$.
\end{lemma}
\begin{proof}
Set $G =  \PSL(2,\bF_5)$, and view $V$ as a $G$-module. Let $U\le G$ be the image of the subgroup of strictly upper triangular matrices. Lemmas~\ref{l.minimalAnIrreducible} and \ref{l.Alt5} show that $V$ is dc-irreducible.

\begin{claim}
For some index set $I$, $V= \bigoplus_{i\in I} V_i$ with each $V_i \cong \Ad(2,\bF_5)$.
\end{claim}
\begin{proofclaim}
We aim to verify the hypotheses of Fact~\ref{fact.Ad(PSL2)}. By  Corollary~\ref{cor.Ad(PSL2)}, we may ignore that $V$ be square free, and 
by dc-irreducibility together with Lemma~\ref{l.C_V(G)Trivial}, we already know $[G,V] = V$ and $C_V(G) = 0$. Further, $U$ is cyclic of prime order, so we certainly have $C_V(u) = C_V(U)$ for all nontrivial $u\in U$. 

It remains to show that $U$ acts cubically on $V$.  As in the proof of Lemma~\ref{l.Alt5Recognition}, we may repeatedly use Lemma~\ref{l.primality} together with the fact that $\dim V = 3$ to see that $[U,U,U,V] = 0$, and since the action is faithful, $[U,V] \neq 0$. Now, if $[U,U,V] = 0$, then $U$ acts quadratically on $V$, but then Fact~\ref{fact.Nat(SL2)} implies that $G \cong \SL(2,\bF_5)$, a contradiction. We conclude that $U$ acts cubically.
\end{proofclaim}

\begin{claim}\label{cl.Ad(2,L)}
$V\cong \Ad(2, L)$ for $L = [U,U,V]= C_V(U)$.
\end{claim}
\begin{proofclaim}
We may assume $V_i = \Ad(2,\bF_5)$ in the decomposition $V = \bigoplus_{i\in I} V_i$. From this, it is readily verified that $[U,U,V] = C_V(U)$, which we denote by $L$.
We may then argue as in Claim~\ref{l.Alt5Recognition:cl:Nat(2,L)} of Lemma~\ref{l.Alt5Recognition} to obtain $V \cong \Ad(2, L)$.
%
\end{proofclaim}
This completes the proof.
\end{proof}\setcounter{claim}{0}

\subsection{Minimal modules for \texorpdfstring{$\Alt(6)$}{Alt(6)} in characteristic 3}

Following our approach for $\Mod(\Alt(5),3,5)$ above, we now  identify $V\in \Mod(\Alt(6),3,3)$.

\begin{lemma}
    Let $V\in \Mod(\Alt(6),3,3)$ be faithful. Then $V$ is an $\bF_9$-vector space, and viewing $\Alt(6)$ as $\PSL(2,\bF_9)$, we have $V\cong \Ad(2, L)$ for some $1$-dimensional subgroup $L\le V$.
\end{lemma}
\begin{proof}
Our proof is similar to that of Lemma~\ref{l.Alt5RecognitionDim3Char5}. We will only show that  $V= \bigoplus_{i\in I}  \Ad(\PSL(2,\bF_9))$ for some index set $I$; the proof that $V\cong \Ad(2, L)$ then follows using the same proof as in Claim~\ref{cl.Ad(2,L)} of Lemma~\ref{l.Alt5RecognitionDim3Char5}.

As before, set $G =  \PSL(2,\bF_9)$, view $V$ as a $G$-module, and let $U\le G$ be the image of the subgroup of strictly upper triangular matrices. By Lemmas~\ref{l.minimalAnIrreducible} and \ref{l.Alt6}, $V$ is dc-irreducible.

All conditions of Fact~\ref{fact.Ad(PSL2)} are easily verified as in the proof of Lemma~\ref{l.Alt5RecognitionDim3Char5} with the exception of \eqref{fact.Ad(PSL2).i:C_V(u)},  
so we now work to show   $C_V(u) = C_V(U)$ for  nontrivial $u\in U$. The subgroup $U$ is a Sylow $3$-subgroup of $G$, so we may choose an isomorphism of $G$ with $\Alt(6)$ that maps $U$ to $\langle (123),(456)\rangle$. Identifying $G$ with $\Alt(6)$ in this way, consider the subgroup \[S=\left\{ \rho\in\Alt(6)\mid \rho\left(\{5,6\}\right) = \{5,6\}\right\},\] which is isomorphic to $\Sym(4)$.

We first show that $C_V(\gamma)$ is dim-connected and of dimension $1$ for every $3$-cycle $\gamma$. 
By Lemma~\ref{l.Sym4Recognition}, $V$ is isomorphic as an $S$-module to either $\std(4,L)$ or $\std^\sigma(4,L)$ for some $1$-dimensional $L\le V$. Computing in $\perm(4,\bZ) = \bZ{e_1} \oplus \cdots \oplus \bZ{e_4}$ and using that $V$ has characteristic $3$, we find that \[\ad_{(123)}\circ\ad_{(123)}(V) = \langle e_1 + e_2 + e_3\rangle \otimes L =  C_V((123));\] thus, (using Fact~\ref{fact.ConnectednessProps}) we see that $C_V((123))$ is dim-connected and of dimension $1$. As any three cycle $\gamma$ is conjugate to $(123)$, $C_V(\gamma)$ is also dim-connected of dimension $1$. 

We next show that $C_V((123)) = C_V((456))$. Consider here the group \[T=\left\{ \rho\in\Alt(6)\mid \ \rho\left(\{1,2\}\right) = \{1,2\}, \rho(3) = 3\right\},\]
which contains $(456)$ and is isomorphic to $\Sym(3)$. Every element of $T$ either centralizes or inverts $(123)$, so $T$ acts on $C_V((123))$. Since $\dim C_V((123)) = 1$, Lemma~\ref{l.Sym3} shows that this action is \emph{not} faithful, implying that $(456)$ acts trivially on $C_V((123))$. Thus, $C_V((123)) \le C_V((456))$, so as both centralizers are dim-connected of dimension $1$, they must in fact be equal.

We currently have that $C_V((123)) = C_V((456)) = C_V(\langle (123),(456)\rangle)$, so it remains to show that if $\varepsilon \in \langle (123),(456)\rangle$ is a product of  disjoint $3$-cycles, then we also have $C_V(\varepsilon) = C_V(\langle (123),(456)\rangle)$. Certainly the former centralizer contains the latter, so it suffices to show that $C_V(\varepsilon)$ is dim-connected of dimension $1$. This is quickly seen using that $(123)$ is  conjugate to $\varepsilon$ via an outer automorphism $\alpha$ of $G$ (see \cite{JR82}): replacing $S$ by its image of $S$ under $\alpha$ in our analysis above shows that $C_V(\varepsilon)$ is indeed dim-connected of dimension $1$.
\end{proof}

\section*{Acknowledgements}
Work on this project began in earnest during 2021 while the first and last authors were master's  students at California State University, Sacramento. Theorem~\ref{t.AltBounds} was the initial goal. Theorem~\ref{t.AltRec} had been suggested by the second author for some time, and with delay, it was picked up and realised in 2023. 

The work of the first and third authors were partially supported by the National Science Foundation under grant No.~DMS-1954127. The first author would like to thank his wife and life partner Alana for making a world with him.

\bibliographystyle{alpha}
\bibliography{SymAlt}
\end{document}